\documentclass{elsarticle}

\usepackage{lineno,hyperref}
\usepackage{amsmath}
\usepackage{amsfonts,amssymb}
\usepackage{esint}
\usepackage{subfigure}
\usepackage{wasysym}
\usepackage{amsthm}
\usepackage{multirow}
\usepackage{algorithm}
\usepackage{algorithmic}
\usepackage{listings}

\newtheorem{remark}{Remark}

\modulolinenumbers[5]

\begin{document}

\begin{frontmatter}

\title{A class of ENO schemes with adaptive order for solving hyperbolic conservation laws}


\author[add1]{Hua Shen\corref{mycorrespondingauthor}}
\cortext[mycorrespondingauthor]{Corresponding author}
\ead{huashen@uestc.edu.cn}
\address[add1]{School of Mathematical Sciences, University of Electronic Science and Technology of China, Chengdu, Sichuan 611731, China}

\begin{abstract}
We propose a class of essentially non-oscillatory schemes with adaptive order (ENO-AO) 
for solving hyperbolic conservation laws.
The new schemes select candidate stencils by novel smoothness indicators
which are the measurements of the minimum discrepancy between
the reconstructed polynomials and the neighboring cell averages.
The new smoothness indicators measure the smoothness of candidate stencils with unequal sizes in a uniform way,
so that we can directly use them to select the optimal stencil from candidates that
range from first-order all the way up to the designed high-order.
Some benchmark test cases are carried out to demonstrate the accuracy and robustness of the proposed schemes.

\end{abstract}


\end{frontmatter}


\section{Introduction}\label{sec:intro}
Many physical processes are governed by nonlinear hyperbolic conservation laws of which
the solutions commonly contain both discontinuities and sophisticated structures with multi-scales.
Therefore, numerical schemes with high-order accuracy and excellent shock-capturing properties are desired 
for solving hyperbolic conservation laws.
The essentially non-oscillatory (ENO) schemes and weighted ENO (WENO) schemes are representatives of such schemes
and achieve a great success in many applications.

Harten and his coworkers \cite{Harten1987ENO} originally proposed a class of ENO schemes 
that are uniformly high-order accurate in smooth regions and are oscillation-free near discontinuities. 
The core idea of ENO schemes is to adaptively choose the smoothest stencil from several local candidates.
Shu and Osher \cite{Shu1988EfficientENO, Shu1989EfficientENOII} efficiently implemented ENO schemes
by using Runge-Kutta time discretizations and numerical fluxes reconstructions instead of cell average reconstructions.
Inspired by ENO schemes, Liu \emph{et al}. \cite{Liu1994WENO} proposed WENO schemes 
by assigning proper weights to the local candidate stencils. 
The principle of designing the weights is to ensure the weights for the stencils containing discontinuities are close to zero,
so that the scheme can achieve non-oscillatory feature near discontinuities.
Jiang and Shu \cite{Jiang_Shu1996WENO} provided a general framework for calculating smooth indicators and weights for 
candidate stencils so that the constructed WENO schemes can recover the optimal order in smooth regions 
while maintaining oscillation-free feature at discontinuities.
Henrick \emph{et al}. \cite{Henrick2005WENO_M} realized that the WENO schemes of Jiang and Shu lose accuracy at critical points
and proposed a simple mapping method to cure this issue.
Borges \emph{et al}. \cite{Borges2008WENO_Z} further improved the accuracy of WENO schemes 
by introducing a parameter to measure the global smoothness and redesigning the nonlinear weights.
Up to now, there are a lot of extended ENO/WENO schemes 
and other shock-capturing high-order schemes that are constructed in the light of ENO/WENO idea,
such as the monotonicity preserving WENO schemes \cite{Balsara2000Monotonicity_WENO}, 
the weighted compact nonlinear schemes \cite{Deng2000WCNS},
the central WENO schemes \cite{Bianco1999CENO, Levy1999CWENO}, 
the Hermite WENO schemes \cite{Qiu2004HermiteWENO, Qiu2005HermiteWENO2D}, 
the  WENO-ADER schemes \cite{Toro2002ADER, Titarev2004WENO_FV, Balsara2013ADER}, 
the $P_NP_M$ schemes \cite{Dumbser2008PNPM, Dumbser2009PNPM, Dumbser2010PNPM},
the adaptive central-upwind WENO schemes \cite{Hu2010WENO_CU},
the targeted ENO (TENO) schemes \cite{Fu2016TENO,Fu2017TENO2}, and so forth.

Almost all the above mentioned ENO/WENO schemes are constructed based on candidate stencils with an equal size.
This kind of schemes cannot achieve the optimal accuracy in some cases.
For example, a five-point global stencil contains three three-point sub-stencils, 
and ENO/WENO schemes constructed on equal-size sub-stencils can only utilize the information 
of the smooth stencils among the three candidates.
However, in some cases, a non-smooth five-point global stencil may contain a smooth four-point sub-stencil,
and in some other cases, all the three three-point sub-stencils may be non-smooth.
Therefore, it is useful to construct numerical schemes based on candidate stencils with unequal sizes.
Levy \emph{et al}. \cite{Levy2000CompactCWENO} proposed a compact central WENO reconstruction
which is a nonlinear combination of candidate stencils with unequal sizes \cite{Levy2000CompactCWENO}.
Recently, Zhu \emph{el al.} \cite{Zhu2016NewWENOFD, Zhu2018NewMultiResolutionWENO, Zhu2019NewMultiResolutionWENO_Tri, Zhu2020NewMultiResolutionWENO_Tet} 
and Balsara \emph{et al.} \cite{Balsara2016WENO_AO, Balsara2020WENO_AO_Unstructure}
extensively implemented the technique proposed by Levy \emph{et al}. \cite{Levy2000CompactCWENO}
to construct WENO schemes with adaptive order and multi-resolution.
Realizing the WENO reconstruction proposed by Levy \emph{et al}. \cite{Levy2000CompactCWENO}
does not always guarantee convexity, Shen \cite{Shen2021WENO_AOA} fixed this issue by a simple modification
and applied it to construct weighted compact central schemes (WCCS) \cite{Shen2021WCCS},.
Fu \emph{et al.} \cite{Fu2018TENO_AO} also constructed TENO schemes with adaptive order based on candidate stencils with unequal sizes.

To the best knowledge of the author,
original ENO schemes select condidate stencils by the Newton divided differences,
and almost all the other ENO/WENO schemes are constructed based on the same framework
in the sense that the smoothness indicators for sub-stencils are the measurement of 
the summation of the square of all space derivatives. 
The above two approaches are not suitable to directly compare the smoothness of stencils with unequal sizes.
Although some WENO/ENO schemes may define different smoothness indicators, the fundamental idea keeps the same.
For the ENO/WENO schemes constructed on candidate stencils with unequal sizes
\cite{Levy2000CompactCWENO, Zhu2016NewWENOFD, Zhu2018NewMultiResolutionWENO, Zhu2019NewMultiResolutionWENO_Tri, 
Zhu2020NewMultiResolutionWENO_Tet, Balsara2016WENO_AO, Balsara2020WENO_AO_Unstructure, Fu2018TENO_AO},
in order to avoid losing accuracy, the nonlinear weights should be designed in special ways based on 
the conventional smoothness indicators.
In this paper, we propose a class of essentially non-oscillatory schemes 
with adaptive order (ENO-AO) based on novel smoothness indicators
which are the measurements of the minimum discrepancy between
the reconstructed polynomials and the neighboring cell averages.
The new smoothness indicators measure the smoothness of candidate stencils with unequal sizes in a uniform way,
so that we can directly use them to select the optimal stencil from candidates that
range from first-order all the way up to the designed high-order. 
As a result, the proposed schemes are high-order accurate in smooth regions 
and are oscillation-free near discontinuities.

\section{A short description of finite difference ENO/WENO schemes}\label{SEC:ENO_WENO_FD}
We consider the one-dimensional scalar conservation law,
\begin{equation}\label{Eq:1DHCL_Eq}
    \frac{\partial u}{\partial t}+\frac{\partial f(u)}{\partial x}=0, t\in[0,\infty),
\end{equation}
in the spatial domain $[x_L,x_R]$ that is discretized into uniform intervals
by $x_j=x_L+(j-1)\Delta x$ ($j=1 \text{ to } N+1$), where $\Delta x=(x_R-x_L)/N$.
A conservative finite difference scheme can be expressed as
\begin{equation}\label{Eq:1DWENOFD}
  \frac{du_j(t)}{dt}=\mathcal{L}(u_j(t))=-\frac{\hat{f}_{j+1/2}-\hat{f}_{j-1/2}}{\Delta x},
\end{equation}
where the numerical flux $\hat{f}_{j\pm1/2}$ is an approximation of the function $h(x)$,
that is implicitly defined by
$f(u(x))=\frac{1}{\Delta x}\int_{x-\Delta x/2}^{x+\Delta x/2}h(\xi)d\xi$ \cite{Shu1988EfficientENO}, at $x_{j\pm1/2}$.
In this manner, we can get a high-order approximation of $\frac{\partial f(u)}{\partial x}$ 
by approximating  $\hat{f}_{j\pm1/2}$ to high-order, thereby getting a semi-discretized scheme with high-order space accuracy.
The attractive feature of this approach is the straightforward extension to multi-dimensional cases,
because the same procedure for approximating the numerical fluxes
can be implemented in a dimension-by-dimension manner.
Once we get the semi-discretized scheme,
we use the third-order TVD Runge-Kutta method \cite{Shu1988EfficientENO} to solve the system of ordinary differential equations,
i.e., $\frac{du_j(t)}{dt}=\mathcal{L}(u_j(t))$.

In order to construct a robust scheme,
we usually split the flux into two parts as
\begin{equation}\label{Eq:Flux_Split}
  f(u)=f^+(u)+f^-(u),
\end{equation}
where $\frac{df^+(u)}{du}\ge0$ and $\frac{df^-(u)}{du}\le0$.
Here, we adopt the global Lax–Friedrichs flux splitting method which is expressed as
\begin{equation}\label{Eq:LF_Flux_Split}
  f^\pm(u)=\frac{1}{2}(f(u)\pm\alpha u),
\end{equation}
where $\alpha=\mbox{max}\left|\frac{df(u)}{du}\right|$
and the maximum is taken over the whole computational mesh points.
After splitting the flux, we respectively construct polynomials to approximate $\hat{f}^\pm(u)$ on 
several sub-stencils of an upwind-biased (or a central) large stencil,
and then choose the smoothest sub-stencil (ENO) or 
use a set of non-linear weights to convexly combine all the sub-stencils (WENO).
The property of ENO/WENO schemes heavily relies on the smoothness indicators.
Jiang and Shu \cite{Jiang_Shu1996WENO} proposed a standard way to calculate the smoothness indicators 
of a $n$th order polynomial $p_n(x)$ on $[x_{i-1/2},x_{i+1/2}]$ as
\begin{equation}\label{Eq:Smoothness_Indicator}
  IS_n=\sum_{l=1}^{n}\Delta x^{2l-1}\int_{x_{i-1/2}}^{x_{i+1/2}}\left(\frac{d^lp_n(x)}{dx^l}\right)^2dx,
\end{equation}
which is the cornerstone of ENO/WENO schemes nowadays.
Balsara \emph{et al.} \cite{Balsara2016WENO_AO} provided an efficient way
to construct a polynomial and to calculate the corresponding smoothness indicator 
defined by Eq. (\ref{Eq:Smoothness_Indicator}) on a given stencil by using Legendre basis.
Since the detail of standard ENO/WENO reconstructions can be found in many references,
we omit it here.

\section{ENO schemes with adaptive order}
\subsection{A novel definition of smoothness indicators}
As we can see, the smoothness indicators defined by Eq. (\ref{Eq:Smoothness_Indicator})
measure the summation of the square of all space derivatives, so they are not suitable to directly compare the smoothness
of polynomials with unequal degrees which are constructed on stencils with unequal sizes.
For example, in smooth region, we prefer to select high-order polynomials,
but Eq. (\ref{Eq:Smoothness_Indicator}) tells us that low-order polynomials are smoother
than high-order polynomials under normal circumstances.
According to the idea of ENO, we should choose the smoothest stencil that is a low-order one in smooth regions.
Therefore, if we want to construct high-order ENO schemes by using candidate stencils with unequal sizes,
we must define new smoothness indicators that can measure 
the smoothness of candidate stencils with unequal sizes in a uniform way.

We denote the interval $[x_{j-1/2},x_{j+1/2}]$ as $I_j$ 
and the corresponding cell average as $f_j$. 
Assume that we have construct a (m+n)th-order polynomial $P_{j-m}^{j+n}(x)$ 
on the stencil $S_{j-m}^{j+n}=[I_{j-m},...,I_j,...,I_{j+n}]$ ($m\ge 0,n\ge 0$) to approximate $h(x)$ on $I_j$,
and the $P_{j-m}^{j+n}(x)$ satisfies
\begin{equation}\label{Eq:Polynomial_Construction}
  \int_{I_k} P_{j-m}^{j+n}(x)dx=f_k, \quad k\in[j-m,j+n].
\end{equation}
Then we define the smoothness indicator for $P_{j-m}^{j+n}(x)$ as
\begin{subequations}\label{Eq:New_Smoothness_Indicator}
\begin{equation}
    IS_{j-m}^{j+n}={\rm MIN}(\delta_L,\delta_R),
\end{equation}
\text{where,}
\begin{equation}
  \delta_L=\begin{cases}
    \frac{1}{2}\left(\left|f_j-f_{j-1}\right|+\left|f_{j-1}-f_{j-2}\right|\right), \quad m=n=0,\\
    \left|\frac{1}{\Delta x}\int_{I_{j-m-1}} P_{j-m}^{j+n}(x)dx-f_{j-m-1}\right|, \quad {\rm otherwise},
  \end{cases}
\end{equation}
\begin{equation}
  \delta_R=\begin{cases}
    \frac{1}{2}\left(\left|f_j-f_{j+1}\right|+\left|f_{j+1}-f_{j+2}\right|\right), \quad m=n=0,\\
    \left|\frac{1}{\Delta x}\int_{I_{j+n+1}} P_{j-m}^{j+n}(x)dx-f_{j+n+1}\right|, \quad {\rm otherwise}.
  \end{cases}
\end{equation}
\end{subequations}

The new smoothness indicator is a measurement of the minimum discrepancy between $P_{j-m}^{j+n}(x)$ 
and the cell averages at the neighbors of $S_{j-m}^{j+n}$.
It has the following properties:
\begin{itemize}
  \item Since $P_{j-m}^{j+n}(x)$ is a $(m+n)$th-order approach to $h(x)$, 
  if $h(x)\in C^{m+n+1}(\mathbb{R})$ on $S_{j-m}^{j+n}$, then we have 
  $IS_{j-m}^{j+n}\propto O(\Delta x^{m+n+1})$ which indicates that 
  high-order polynomials are usually smoother than low-order ones.
  \item If $S_{j-m}^{j+n}$ contains any discontinuity, $IS_{j-m}^{j+n}\propto O(1)$ which is larger than the
  smoothness indicators of its smooth sub-stencils.
\end{itemize}
Therefore, if we use the new smoothness indicators to select candidate stencils with unequal sizes,
we can pick up the highest-order candidate in smooth regions and the smoothest candidate near discontinuities.
\begin{remark}
  We note that, the definition of $IS_j^j$  includes the difference between neighbor and neighbor's neighbor
  which is a little bit different from others.
This operation is used to avoid lost of accuracy at first-order extrema.
If we only use the difference between the target cell and its neighbor,
$IS_j^j$ may become zero at first-order extrema.
Then we will mistakenly select $P_j^j$ as the `smoothest' candiate 
which has the lowest-order omong the candidates.
\end{remark}

\subsection{ENO reconstructions by using candidate stencils with unequal sizes}
We provide the reconstruction for $\hat{f}_{j+1/2}^+$ and drop the superscript $+$ for simplicity.
The formulas for $\hat{f}_{j-1/2}^-$ are symmetric with respect to $x_{j+1/2}$.
By using the condition given by Eq. (\ref{Eq:Polynomial_Construction}),
we can trivially construct $P_{j-m}^{j+n}(x)$ on $S_{j-m}^{j+n}$ in the following form
\begin{equation}
  P_{j-m}^{j+n}(x)=\sum_{k=0}^{m+n} a_kL_k(x),
\end{equation}
where $L_k(x)$ are Legendre basis.
The specific forms of $P_{j-m}^{j+n}(x)$ on all the sub-stencils of $S_{j-3}^{j+3}$
can be found in Balsara \emph{et al.} \cite{Balsara2016WENO_AO}
and will be not shown here. 
Once we construct $P_{j-m}^{j+n}(x)$, 
we can immediately get $\hat{f}_{j+1/2}$ as given by Table \ref{Table:Flux_reconstruction}.
For the seventh-order reconstruction on $S_{j-3}^{j+3}$,
there are 16 stencils in total.
Before executing the space reconstruction,
we perform Fourier analysis for the linear flux on every stencil given by Table \ref{Table:Flux_reconstruction}.
The real part and imaginary part of the modified wavenumber \cite{Lele1992CompactFD} 
associated with each flux are given by 
Table \ref{Table:Wavenumber_Real} and Table \ref{Table:Wavenumber_Real} respectively.
Since the new smoothness indicator does not take account of all space derivatives,
it cannot completely suppress the numerical instabilities induced by the linearly unstable stencils 
as the traditional ENO/WENO schemes do.
As a result, we only select linearly stable fluxes 
(the imaginary part of the modified wavenumber is nonpositive for all wavenumbers) as candidates. 
We note that, although the flux constructed on $S_{j-1}^j$ is linearly stable,
we abandon it due to its large dispersion errors.
As a result,
the candidate stencils for the fifth-order reconstruction include 
$S_j^j$, $S_j^{j+1}$, $S_{j-1}^{j+1}$, $S_{j-2}^{j+1}$, $S_{j-1}^{j+2}$, $S_{j-2}^{j+2}$.
The candidate stencils for the seventh-order reconstruction include
$S_j^j$, $S_j^{j+1}$, $S_{j-1}^{j+1}$, $S_{j-2}^{j+1}$, $S_{j-1}^{j+2}$, $S_{j-2}^{j+2}$, 
$S_{j-3}^{j+2}$, $S_{j-2}^{j+3}$, $S_{j-3}^{j+3}$.
The specific forms of $\delta_L$ and $\delta_R$ in Eq. (\ref{Eq:New_Smoothness_Indicator}) 
associated with the candidate stencils
are respectively given by Table \ref{Table:DeltaL} and Table \ref{Table:DeltaR}.

The reconstruction principle is as follows:
\begin{itemize}
  \item If multiple stencils' smoothness indicators are smaller than a small number $\delta$,
  the highest-order stencil has the highest priority, 
  and the central stencil has a higher priority than the upwind-biased stencil with the same order.
  \item Otherwise, we choose the stencil with the smallest smoothness indicator.
\end{itemize}
To make the whole reconstruction procedure clear,
we attach the C code of seventh-order ENO-AO reconstruction in \ref{AppendixA}.

\begin{remark}
  The first reconstruction principle is used to avoid lost of optimal accuracy in smooth regions.
  The parameter $\delta$ has a similar function as the parameter $\varepsilon$ in WENO schemes.
  In WENO schemes, the weights have the form of $w_n=\left(\frac{1}{IS_n+\varepsilon}\right)^p$ \cite{Jiang_Shu1996WENO}.
  Therefore, when $IS_n<<\varepsilon$, $w_n\thickapprox \frac{1}{\varepsilon}$ and
  the effect of smoothness indicators vanishes.
  In WENO-JS \cite{Jiang_Shu1996WENO}, $\varepsilon=10^{-6}$.
  In this paper we set $\delta=10^{-5}$. 
  We note that, the smoothness indicators of WENO schemes Eq. (\ref{Eq:Smoothness_Indicator}) have the dimension of $(\Delta f)^2$,
  and the the new smoothness indicators Eq. (\ref{Eq:New_Smoothness_Indicator}) have the dimension of $\Delta f$.
  When $IS_n\le0.01\varepsilon$, smoothness indicators have a weak effect on the weights of WENO schemes.
  Therefore, $\delta=10^{-5}$ is approximate to $\varepsilon=10^{-8}$ in WENO schemes
  which is small enough under normal conditions.

\end{remark}

\begin{table}
  \centering
  \begin{tabular}{|c|l|}
     \hline
                    & $\hat{f}_{j+1/2}$     \\
     \hline
     $S_j^j$        &$f_j$\\
     $S_j^{j+1}$     &$(f_j+f_{j+1})/2$\\
     $S_j^{j+2}$     &$(2f_j+5f_{j+1}-f_{j+2})/6$\\
     $S_j^{j+3}$     &$(3f_j+13f_{j+1}-5f_{j+2}+f_{j+3})/12$\\
     $S_{j-1}^j$     &$(-f_{j-1}+3f_j)/2$\\
     $S_{j-1}^{j+1}$     &$(-f_{j-1}+5f_j+2f_{j+1})/6$\\
     $S_{j-1}^{j+2}$     &$(-f_{j-1}+7f_j+7f_{j+1}-f_{j+2})/12$\\
     $S_{j-1}^{j+3}$     &$(-3f_{j-1}+27f_j+47f_{j+1}-13f_{j+2}+2f_{j+3})/60$\\
     $S_{j-2}^{j}$     &$(2f_{j-2}-7f_{j-1}+11f_j)/6$\\
     $S_{j-2}^{j+1}$     &$(f_{j-2}-5f_{j-1}+13f_j+3f_{j+1})/12$\\
     $S_{j-2}^{j+2}$     &$(2f_{j-2}-13f_{j-1}+47f_j+27f_{j+1}-3f_{j+2})/60$\\
     $S_{j-2}^{j+3}$     &$(f_{j-2}-8f_{j-1}+37f_j+37f_{j+1}-8f_{j+2}+f_{j+3})/60$\\
     $S_{j-3}^{j}$     &$(-3f_{j-3}+13f_{j-2}-23f_{j-1}+25f_j)/12$\\
     $S_{j-3}^{j+1}$     &$(-3f_{j-3}+17f_{j-2}-43f_{j-1}+77f_j+12f_{j+1})/60$\\
     $S_{j-3}^{j+2}$     &$(-f_{j-3}+7f_{j-2}-23f_{j-1}+57f_j+22f_{j+1}-2f_{j+2})/60$\\
     $S_{j-3}^{j+3}$     &$(-3f_{j-3}+25f_{j-2}-101f_{j-1}+319f_j+214f_{j+1}-38f_{j+2}+4f_{j+3})/420$\\

     \hline
  \end{tabular}
  \caption{The numerical flux $\hat{f}_{j+1/2}$ on different stencils.}
  \label{Table:Flux_reconstruction}
\end{table}

\begin{table}
  \centering
  \begin{tabular}{|c|l|}
     \hline
                    & Real part of the modified wavenumber $\kappa'$    \\
     \hline
     $S_j^j$        &$sin\kappa$\\
     $S_j^{j+1}$     &$sin\kappa$\\
     $S_j^{j+2}$     &$\frac{4}{3}sin\kappa-\frac{1}{6}sin(2\kappa)$\\
     $S_j^{j+3}$     &$\frac{7}{4}sin\kappa-\frac{1}{2}sin(2\kappa)+\frac{1}{12}sin(3\kappa)$\\
     $S_{j-1}^j$     &$2sin\kappa-\frac{1}{2}sin(2\kappa)$\\
     $S_{j-1}^{j+1}$     &$\frac{4}{3}sin\kappa-\frac{1}{6}sin(2\kappa)$\\
     $S_{j-1}^{j+2}$     &$\frac{4}{3}sin\kappa-\frac{1}{6}sin(2\kappa)$\\
     $S_{j-1}^{j+3}$     &$\frac{3}{2}sin\kappa-\frac{3}{10}sin(2\kappa)+\frac{1}{30}sin(3\kappa)$\\
     $S_{j-2}^{j}$     &$3sin\kappa-\frac{3}{2}sin(2\kappa)+\frac{1}{3}sin(3\kappa)$\\
     $S_{j-2}^{j+1}$     &$\frac{7}{4}sin\kappa-\frac{1}{2}sin(2\kappa)+\frac{1}{12}sin(3\kappa)$\\
     $S_{j-2}^{j+2}$     &$\frac{3}{2}sin\kappa-\frac{3}{10}sin(2\kappa)+\frac{1}{30}sin(3\kappa)$\\
     $S_{j-2}^{j+3}$     &$\frac{3}{2}sin\kappa-\frac{3}{10}sin(2\kappa)+\frac{1}{30}sin(3\kappa)$\\
     $S_{j-3}^{j}$     &$4sin\kappa-3sin(2\kappa)+\frac{4}{3}sin(3\kappa)-\frac{1}{4}sin(4\kappa)$\\
     $S_{j-3}^{j+1}$     &$\frac{11}{5}sin\kappa-sin(2\kappa)+\frac{1}{3}sin(3\kappa)-\frac{1}{20}sin(4\kappa)$\\
     $S_{j-3}^{j+2}$     &$\frac{26}{15}sin\kappa-\frac{8}{15}sin(2\kappa)+\frac{2}{15}sin(3\kappa)-\frac{1}{60}sin(4\kappa)$\\
     $S_{j-3}^{j+3}$     &$\frac{8}{5}sin\kappa-\frac{2}{5}sin(2\kappa)+\frac{8}{105}sin(3\kappa)-\frac{1}{140}sin(4\kappa)$\\

     \hline
  \end{tabular}
  \caption{The real part the modified wavenumber 
  associated with each flux constructed on $S_{j-m}^{j+n}$. 
  Here, $\kappa=k\Delta x$ is the scaled wavenumber.}
  \label{Table:Wavenumber_Real}
\end{table}

\begin{table}
  \centering
  \begin{tabular}{|c|l|}
     \hline
                    & Imaginary part of the modified wavenumber $\kappa'$  \\
     \hline
     $S_j^j$        &$cos\kappa-1$\\
     $S_j^{j+1}$     &$0$\\
     $S_j^{j+2}$     &$\frac{1}{2}-\frac{2}{3}cos\kappa+\frac{1}{6}cos(2\kappa)$\\
     $S_j^{j+3}$     &$\frac{5}{6}-\frac{5}{4}cos\kappa+\frac{1}{2}cos(2\kappa)-\frac{1}{12}cos(3\kappa)$\\
     $S_{j-1}^j$     &$-\frac{3}{2}+2cos\kappa-\frac{1}{2}cos(2\kappa)$\\
     $S_{j-1}^{j+1}$     &$-\frac{1}{2}+\frac{2}{3}cos\kappa-\frac{1}{6}cos(2\kappa)$\\
     $S_{j-1}^{j+2}$     &$0$\\
     $S_{j-1}^{j+3}$     &$\frac{1}{3}-\frac{1}{2}cos\kappa+\frac{1}{5}cos(2\kappa)-\frac{1}{30}cos(3\kappa)$\\
     $S_{j-2}^{j}$     &$-\frac{11}{6}+3cos\kappa-\frac{3}{2}cos(2\kappa)+\frac{1}{3}cos(3\kappa)$\\
     $S_{j-2}^{j+1}$     &$-\frac{5}{6}+\frac{5}{4}cos\kappa-\frac{1}{2}cos(2\kappa)+\frac{1}{12}cos(3\kappa)$\\
     $S_{j-2}^{j+2}$     &$-\frac{1}{3}+\frac{1}{2}cos\kappa-\frac{1}{5}cos(2\kappa)+\frac{1}{30}cos(3\kappa)$\\
     $S_{j-2}^{j+3}$     &$0$\\
     $S_{j-3}^{j}$     &$-\frac{25}{12}+4cos\kappa-3cos(2\kappa)+\frac{4}{3}cos(3\kappa)-\frac{1}{4}cos(4\kappa)$\\
     $S_{j-3}^{j+1}$     &$-\frac{13}{12}+\frac{9}{5}cos\kappa-cos(2\kappa)+\frac{1}{3}cos(3\kappa)-\frac{1}{20}cos(4\kappa)$\\
     $S_{j-3}^{j+2}$     &$-\frac{7}{12}+\frac{14}{15}cos\kappa-\frac{7}{15}cos(2\kappa)+\frac{2}{15}cos(3\kappa)-\frac{1}{60}cos(4\kappa)$\\
     $S_{j-3}^{j+3}$     &$-\frac{1}{4}+\frac{2}{5}cos\kappa-\frac{1}{5}cos(2\kappa)+\frac{2}{35}cos(3\kappa)-\frac{1}{140}cos(4\kappa)$\\

     \hline
  \end{tabular}
  \caption{The imaginary part the modified wavenumber 
  associated with each flux constructed on $S_{j-m}^{j+n}$.
  Here, $\kappa=k\Delta x$ is the scaled wavenumber.}
  \label{Table:Wavenumber_Im}
\end{table}

\begin{table}
  \centering
  \begin{tabular}{|c|l|}
     \hline
                    & $\delta_L$ in Eq. (\ref{Eq:New_Smoothness_Indicator})     \\
     \hline
     $S_j^j$        &$(|f_j-f_{j-1}|+|f_{j-1}-f_{j-2}|)/2$\\
     $S_j^{j+1}$     &$|f_{j-1}-2f_j+f_{j+1}|$\\
     $S_{j-1}^{j+1}$     &$|-f_{j-2}+3f_{j-1}-3f_j+f_{j+1}|$\\
     $S_{j-2}^{j+1}$     &$|f_{j-3}-4f_{j-2}+6f_{j-1}-4f_j+f_{j+1}|$\\
     $S_{j-1}^{j+2}$     &$|f_{j-2}-4f_{j-1}+6f_j-4f_{j+1}+f_{j+2}|$\\
     $S_{j-2}^{j+2}$     &$|-f_{j-3}+5f_{j-2}-10f_{j-1}+10f_j-5f_{j+1}+f_{j+2}|$\\
     $S_{j-3}^{j+2}$     &$|f_{j-4}-6f_{j-3}+15f_{j-2}-20f_{j-1}+15f_j-6f_{j+1}+f_{j+2}|$\\
     $S_{j-2}^{j+3}$     &$|f_{j-3}-6f_{j-2}+15f_{j-1}-20f_j+15f_{j+1}-6f_{j+2}+f_{j+3}|$\\
     $S_{j-3}^{j+3}$     &$|-f_{j-4}+7f_{j-3}-21f_{j-2}+35f_{j-1}-35f_j+21f_{j+1}-7f_{j+2}+f_{j+3}|$\\

     \hline
  \end{tabular}
  \caption{$\delta_L$ in Eq. (\ref{Eq:New_Smoothness_Indicator}) on different stencils.}
  \label{Table:DeltaL}
\end{table}

\begin{table}
  \centering
  \begin{tabular}{|c|l|}
     \hline
                    & $\delta_R$ in Eq. (\ref{Eq:New_Smoothness_Indicator})    \\
     \hline
     $S_j^j$        &$(|f_j-f_{j+1}|+|f_{j+1}-f_{j+2}|)/2$\\
     $S_j^{j+1}$     &$|f_{j}-2f_{j+1}+f_{j+2}|$\\
     $S_{j-1}^{j+1}$     &$|-f_{j-1}+3f_{j}-3f_{j+1}+f_{j+2}|$\\
     $S_{j-2}^{j+1}$     &$|f_{j-2}-4f_{j-1}+6f_{j}-4f_{j+1}+f_{j+2}|$\\
     $S_{j-1}^{j+2}$     &$|f_{j-1}-4f_{j}+6f_{j+1}-4f_{j+2}+f_{j+3}|$\\
     $S_{j-2}^{j+2}$     &$|-f_{j-2}+5f_{j-1}-10f_{j}+10f_{j+1}-5f_{j+2}+f_{j+3}|$\\
     $S_{j-3}^{j+2}$     &$|f_{j-3}-6f_{j-2}+15f_{j-1}-20f_{j}+15f_{j+1}-6f_{j+2}+f_{j+3}|$\\
     $S_{j-2}^{j+3}$     &$|f_{j-2}-6f_{j-1}+15f_{j}-20f_{j+1}+15f_{j+2}-6f_{j+3}+f_{j+4}|$\\
     $S_{j-3}^{j+3}$     &$|-f_{j-3}+7f_{j-2}-21f_{j-1}+35f_{j}-35f_{j+1}+21f_{j+2}-7f_{j+3}+f_{j+4}|$\\

     \hline
  \end{tabular}
  \caption{$\delta_R$ in Eq. (\ref{Eq:New_Smoothness_Indicator}) on different stencils.}
  \label{Table:DeltaR}
\end{table}

\section{Numerical examples}\label{SEC:_NumExam}
In this section, we use some benchmarks to test the properties of the proposed 5th- and 7th-order ENO-AO schemes.
All results are compared with 5th- and 7th-order WENO-Z schemes \cite{Borges2008WENO_Z, Castro2011High_Order_WENO_Z} 
which are improved versions of WENO-JS schemes \cite{Jiang_Shu1996WENO}.
As suggested by the inventors, we set $\epsilon=10^{-40}$, and $p=1$ in the WENO-Z computations
and use the global optimal order smoothness indicator $\tau_{2r-1}^{opt}$ \cite{Castro2011High_Order_WENO_Z}.
In all simulations, we set $CFL=0.3$ except for the spectral analysis and convergence tests.
For computations of the Euler equations,
the Roe average at the cell face is adopted for characteristic decomposition.

\subsection{Numerical examples for the 1D linear advection equation}
\begin{figure}
  \centering
  \subfigure[Dispersion of 5th-order schemes]{
  \label{FIG:5th_order_dispersion}
  \includegraphics[width=5.5 cm]{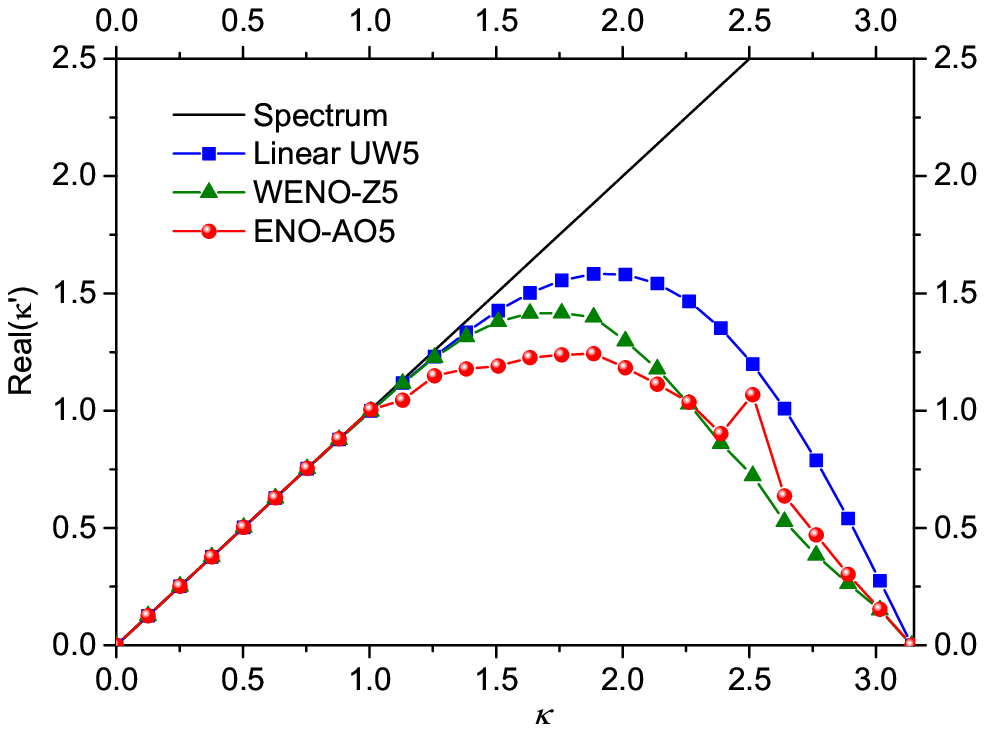}}
  \subfigure[Enlarged view of dispersion errors of 5th-order schemes]{
  \label{FIG:enlarged_5th_order_dispersion}
  \includegraphics[width=5.5 cm]{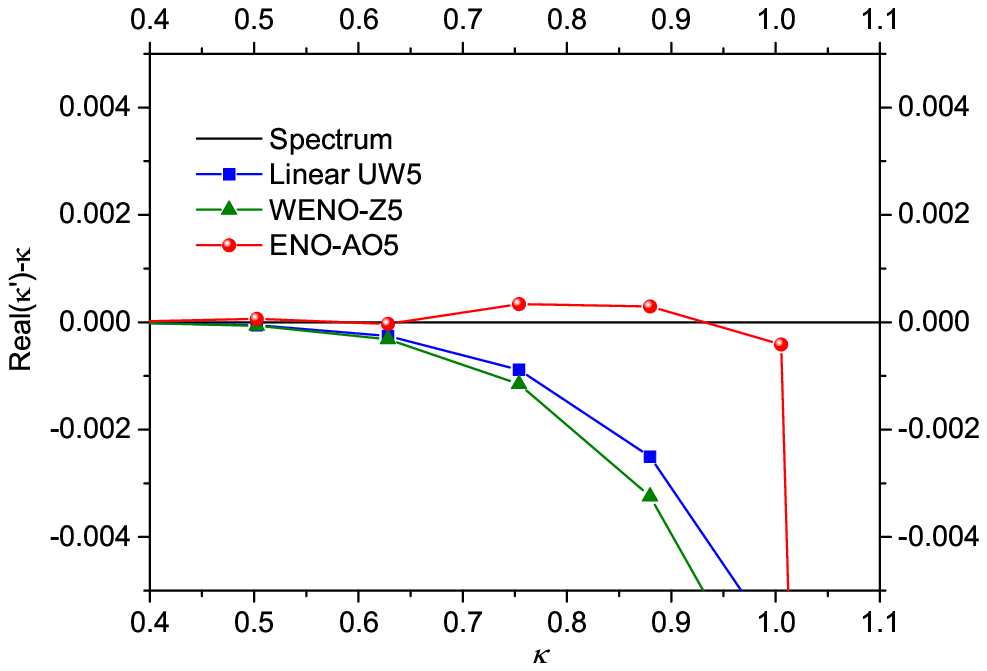}}
  \subfigure[Dissipation of 5th-order schemes]{
  \label{FIG:5th_order_dissipation}
  \includegraphics[width=5.5 cm]{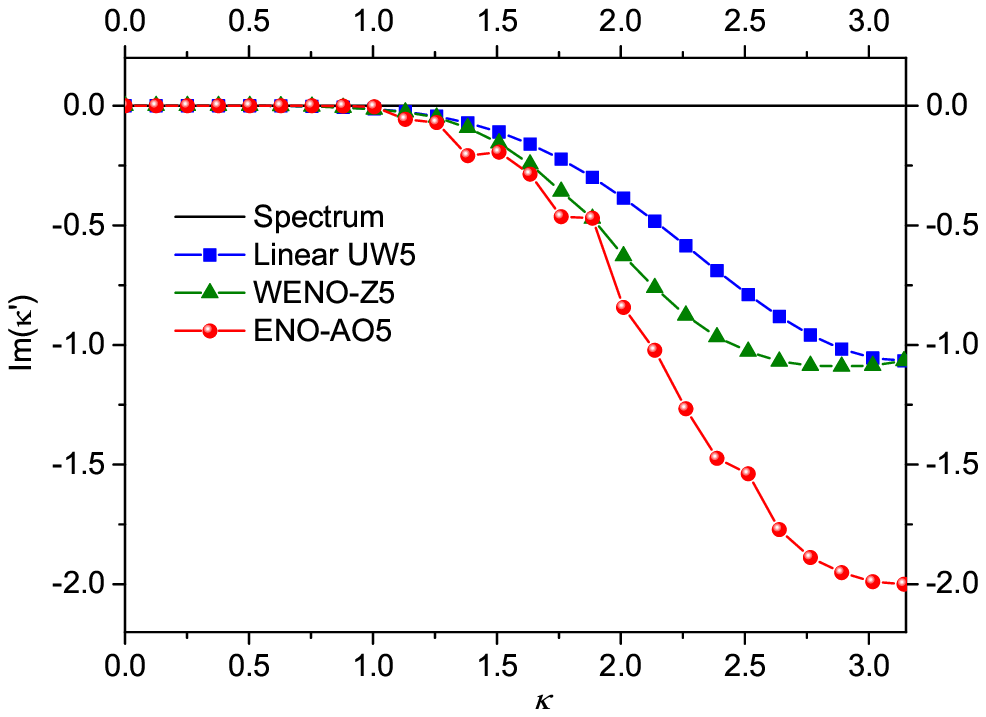}}
  \subfigure[Enlarged view of dissipation of 5th-order schemes]{
  \label{FIG:enlarged_5th_order_dissipation}
  \includegraphics[width=5.5 cm]{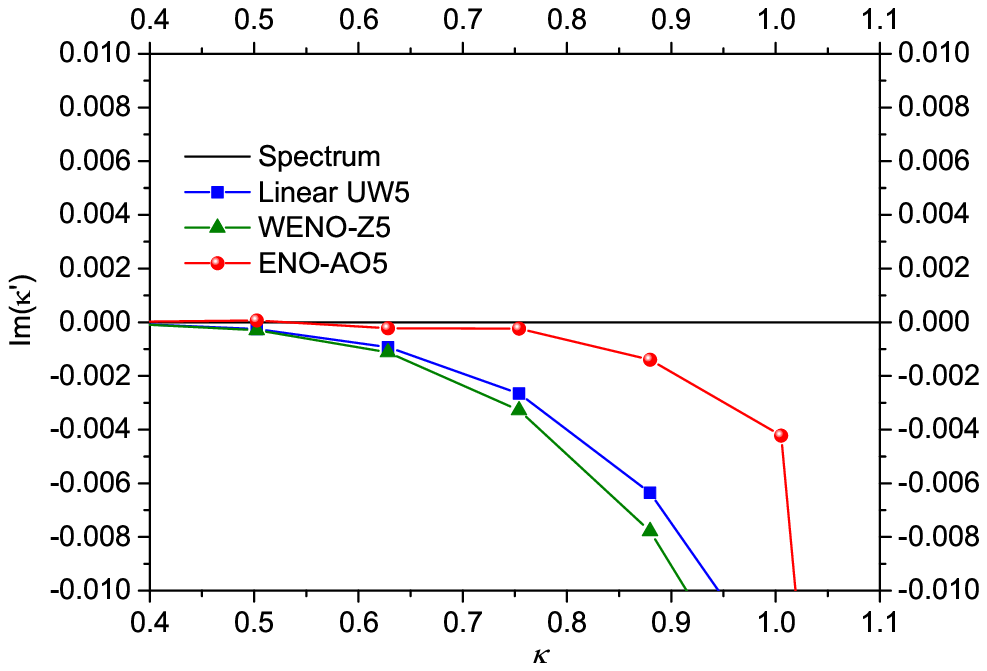}}
  \subfigure[Dispersion of 7th-order schemes]{
  \label{FIG:7th_order_dispersion}
  \includegraphics[width=5.5 cm]{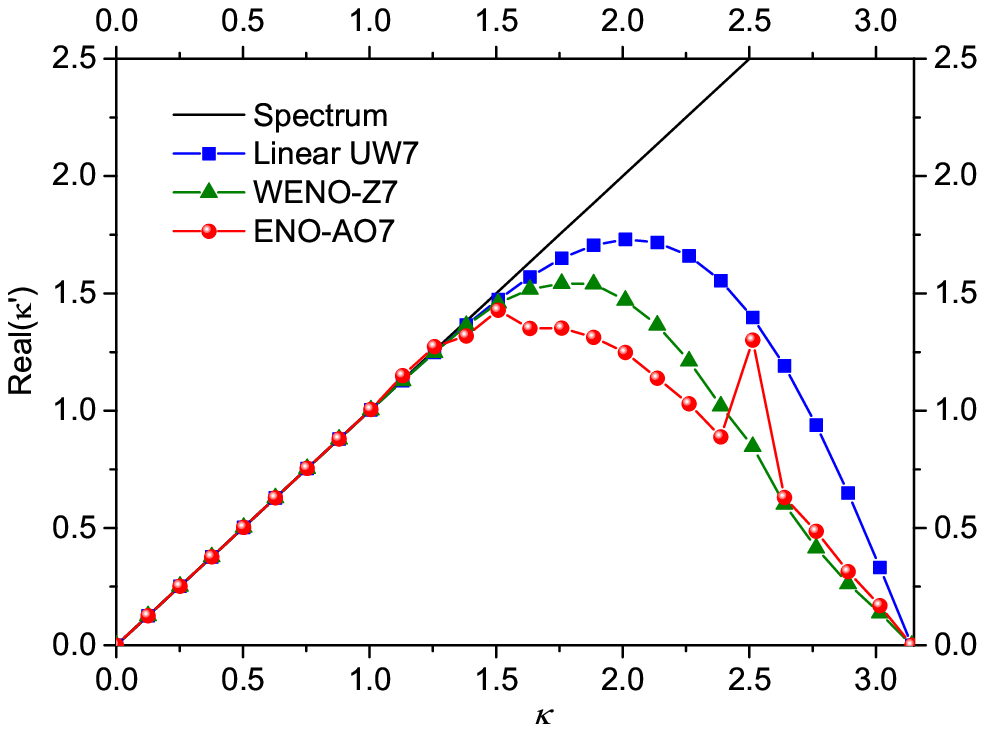}}
  \subfigure[Enlarged view of dispersion errors of 7th-order schemes]{
  \label{FIG:enlarged_7th_order_dispersion}
  \includegraphics[width=5.5 cm]{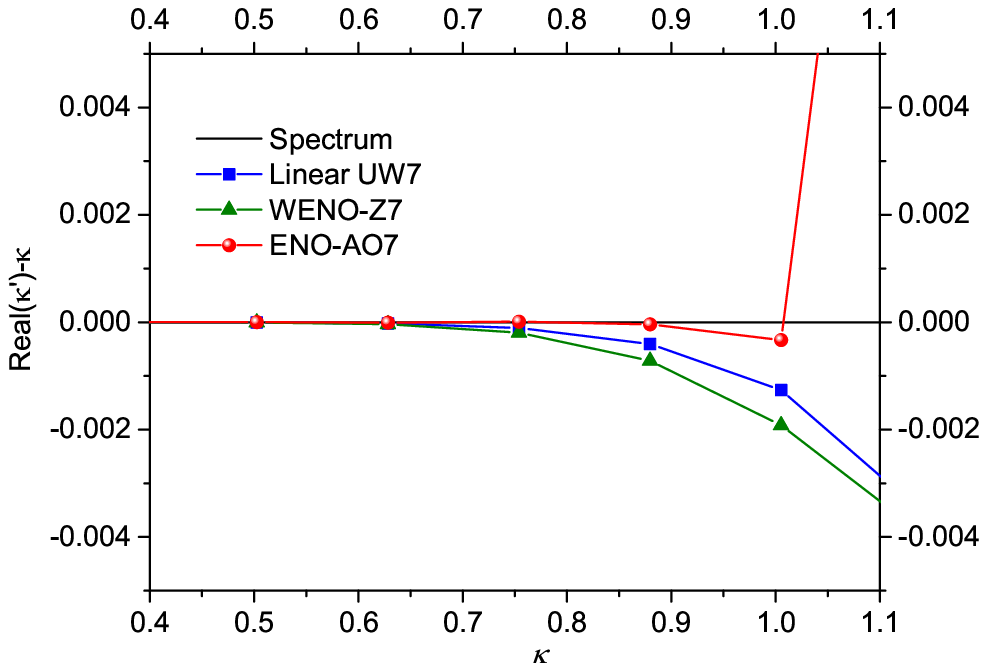}}
  \subfigure[Dissipation of 7th-order schemes]{
  \label{FIG:7th_order_dissipation}
  \includegraphics[width=5.5 cm]{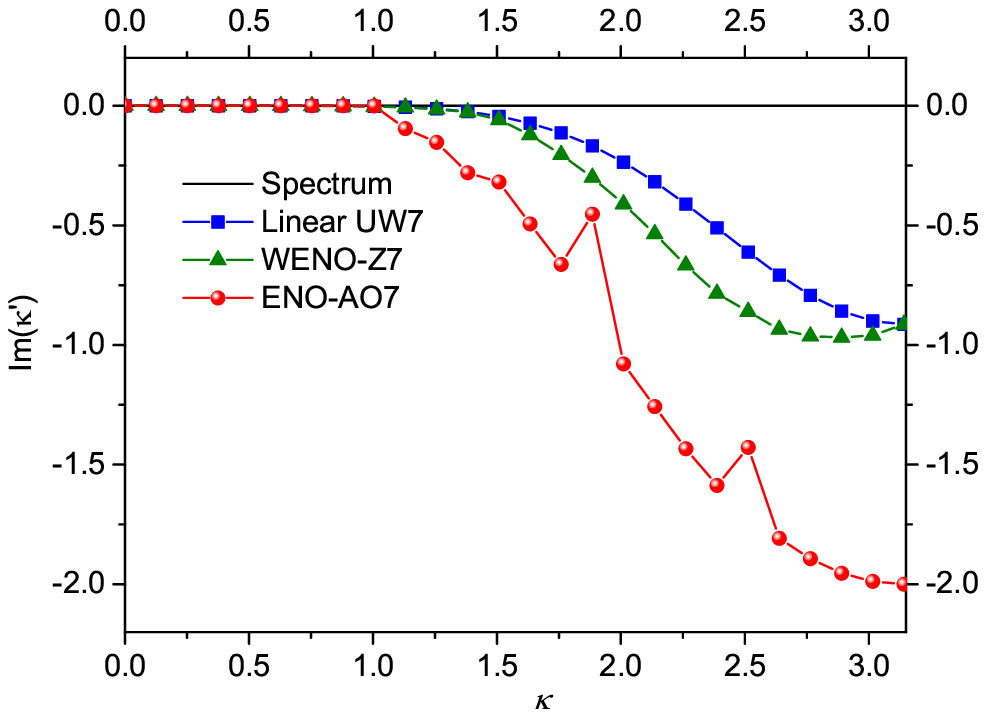}}
  \subfigure[Enlarged view of dissipation of 7th-order schemes]{
  \label{FIG:enlarged_5th_order_dissipation}
  \includegraphics[width=5.5 cm]{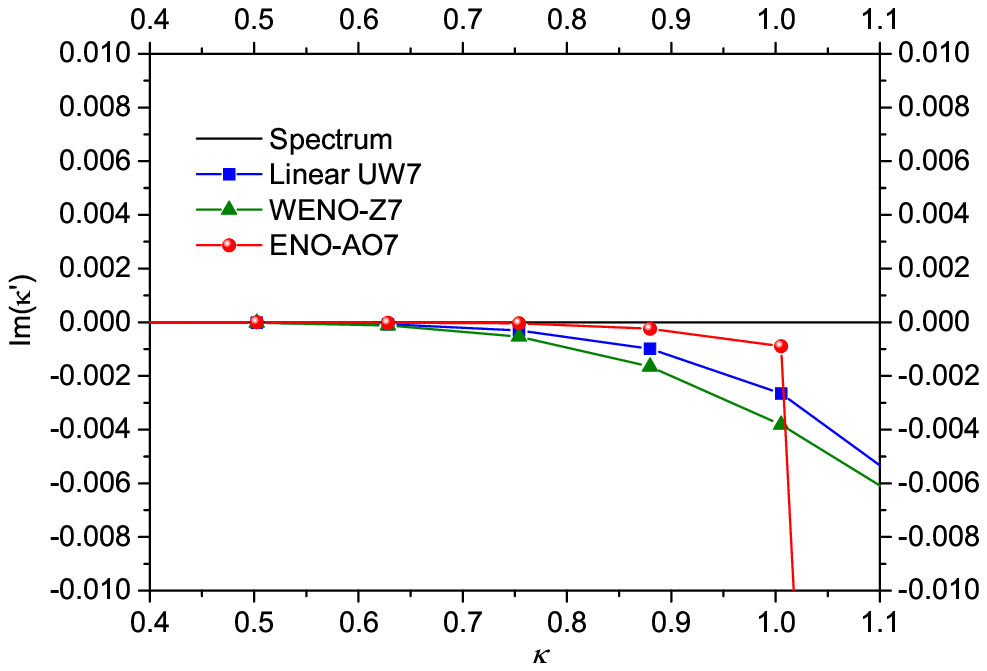}}
  \caption{Comparison of dispersion and dissipation properties of linear upwind (UW), WENO-Z, and ENO-AO schemes.}
\label{FIG:ADR}
\end{figure}

We consider the 1D linear advection equation
\begin{equation}\label{Eq:1DAdvEq}
  \frac{\partial u}{\partial t}+\frac{\partial u}{\partial x}=0.
\end{equation}
We first perform spectrum analysis of different schemes.
Since ENO/WENO schemes are highly nonlinear,
it is difficult to get analytical solutions of modified wavenumbers 
as we do for the linear schemes.
We seek for approximate dispersion relation (ADR) by using a numerical technique
proposed by Pirozzoli \cite{Pirozzoli2006ADR}.
Fig. \ref{FIG:ADR} shows the approximate dispersion and dissipation properties 
of linear upwind (UW), WENO-Z, and ENO-AO schemes.
The real part of modified wavenumber $\kappa'$ represents numerical dispersion,
and the imaginary part represents numerical dissipation.
We note that, in order to easily compare the dispersion properties of different schemes,
we show the dispersion errors $\rm{Real}(\kappa')-\kappa$ 
in Figs. \ref{FIG:enlarged_5th_order_dispersion} and \ref{FIG:enlarged_7th_order_dispersion}.
The left column of Fig. \ref{FIG:ADR} shows that all schemes have very small 
numerical dispersion and dissipation for small wavenumbers
and ENO-AO obviously has larger numerical dispersion and dissipation 
than WENO-Z schemes for large wavenumbers ($\kappa>1$) 
becasue ENO-AO contains lower-order stencils than WENO-Z schemes.
The right column of Fig. \ref{FIG:ADR} shows that  ENO-AO has smaller 
numerical dispersion and dissipation than WENO-Z and linear UW schemes for medium wavenumbers 
($1\gtrsim \kappa\gtrsim0.4$ for 5th-order; $1\gtrsim \kappa\gtrsim0.6$ for 7th-order).
In fact, ($\kappa=1$) is approximately equivalent to 6 elements within a period of a sinusoidal wave
which is considered as a very coarse mesh for most applications.
Therefore, large numerical dissipation for $\kappa>1$ is not a fatal weakness.
On the contrary, it is benificial to eliminate high-frequency oscillations.
Even if we use UW or WENO-Z schemes in the high-wavenumber cases,
we cannot obtain satisfied results.
As we will see in the following, the excellent spectral properties of ENO-AO schemes for the medium wavenumbers
are helpful for resolving small-scale structures in most cases.

Next we test the convergence rate by the advection of a sinusoidal wave. The computational domain is $[-1,1]$,
the initial condition is $u_0(x)=sin(\pi x)$,
and periodic boundary conditions are applied on both sides.
Since we use third-order TVD RK method for time discretization,
the time accuracy does not match with space accuracy.
We respectively set $\Delta t=\Delta x^{5/3}$ and $\Delta t=\Delta x^{7/3}$ for 5th- and 7th-order
schemes to test the convergence rates of spatial operators.
Table \ref{Table:1DAdvectionConvergence} shows the numerical errors of WENO-Z and ENO-AO at $t=2$,
where $L_1$ is the average error and $L_\infty$ is the maximum error.
We observe that ENO-AO achieve the optimal convergence rate for $L_1$ and $L_\infty$
when refining the mesh, and only the errors of ENO-AO5 scheme on coarse meshes are
slightly larger than WENO-Z5 schemes.
It means that the low-order stencils of ENO-AO schemes are not activated for this smooth case.

\begin{table}
  \centering
  \begin{tabular}{|c|c|c|c|c|c|c|c|c|}
    \hline
      & \multicolumn{4}{|c|}{WENO-Z5} & \multicolumn{4}{|c|}{ENO-AO5} \\
     \hline
      Mesh size       & $L_1$  &Order & $L_\infty$  &Order    & $L_1$  &Order & $L_\infty$  &Order    \\
     \hline
     1/20              &7.78E-6 &-     &1.26E-5      &-        &2.47E-5 &-     &6.23E-5      &- \\

     1/30             &1.03E-6 &4.99  &1.66E-6     &5.00     &1.30E-6 &7.26  &3.04E-6      &7.45 \\

     1/40             &2.46E-7 &4.98  &3.93E-7      &5.01     &2.46E-7 &5.79  &3.91E-7      &7.13 \\

     1/50             &8.10E-8 &4.99  &1.29E-7      &4.99     &8.10E-8 &4.99  &1.28E-7      &5.00 \\
     \hline

     & \multicolumn{4}{|c|}{WENO-Z7} & \multicolumn{4}{|c|}{ENO-AO7} \\
     \hline
      Mesh size       & $L_1$  &Order & $L_\infty$  &Order    & $L_1$  &Order & $L_\infty$  &Order    \\
     \hline
     1/20              &3.68E-8 &-     &6.53E-8      &-        &3.68E-8 &-     &5.89E-8      &-  \\

     1/30             &2.17E-9 &6.98  &3.62E-9      &7.13     &2.17E-9 &6.98  &3.46E-9      &6.99 \\

     1/40             &2.92E-10 &6.97  &4.74E-10      &7.07     &2.92E-10 &6.97  &4.63E-10   &6.99 \\

     1/50             &6.19E-11 &6.95  &9.95E-11      &7.00     &6.19E-11 &6.95  &9.81E-11    &6.95 \\
     \hline
  \end{tabular}
  \caption{Numerical errors of the advection of a sinusoidal wave at $t=2$ computed by WENO-Z and ENO-AO with different mesh sizes.}
  \label{Table:1DAdvectionConvergence}
\end{table}

\begin{figure}
  \centering
  \subfigure[WENO-Z5]{
  \label{FIG:Linear_Advection_WENO-Z5}
  \includegraphics[width=5.5 cm]{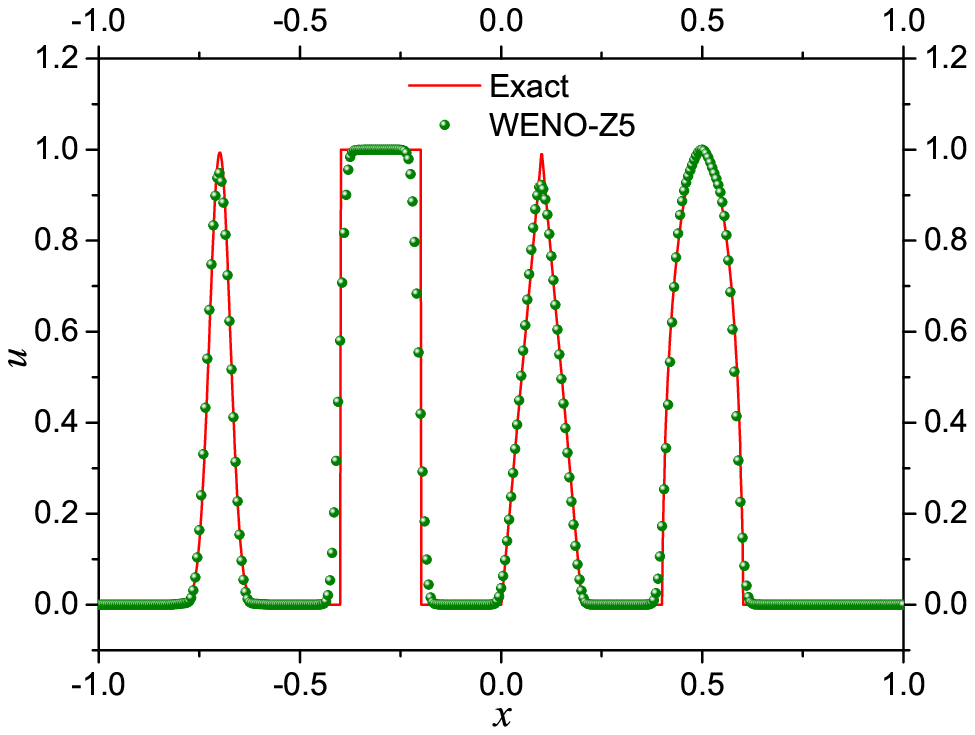}}
  \subfigure[ENO-AO5]{
  \label{FIG:Linear_Advection_ENO-AO5}
  \includegraphics[width=5.5 cm]{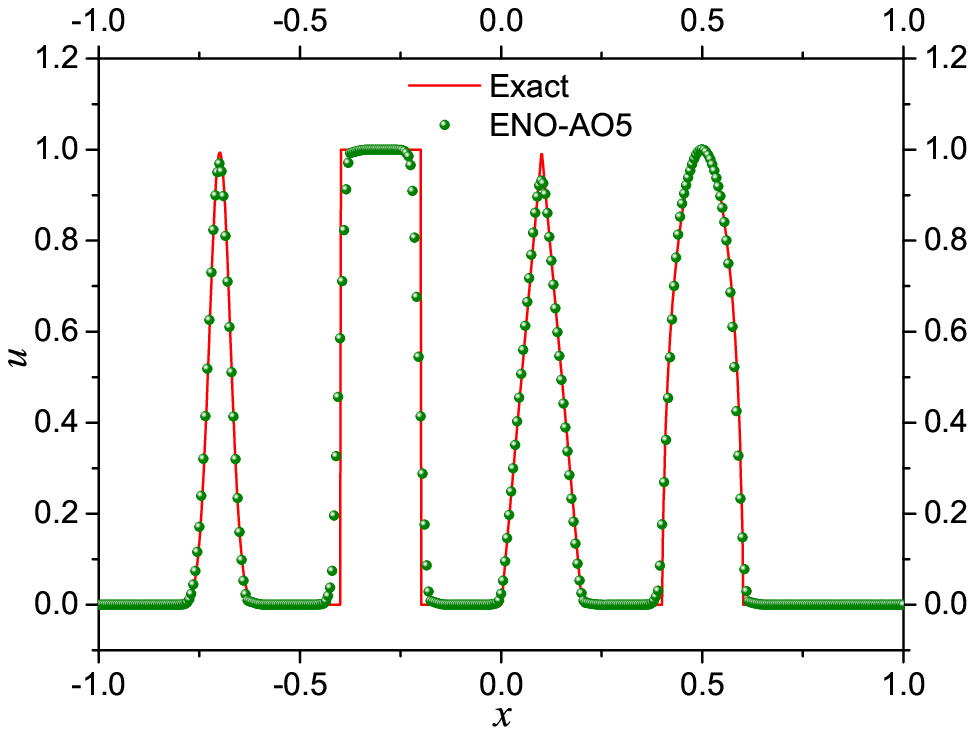}}
  \subfigure[WENO-Z7]{
  \label{FIG:Linear_Advection_WENO-Z7}
  \includegraphics[width=5.5 cm]{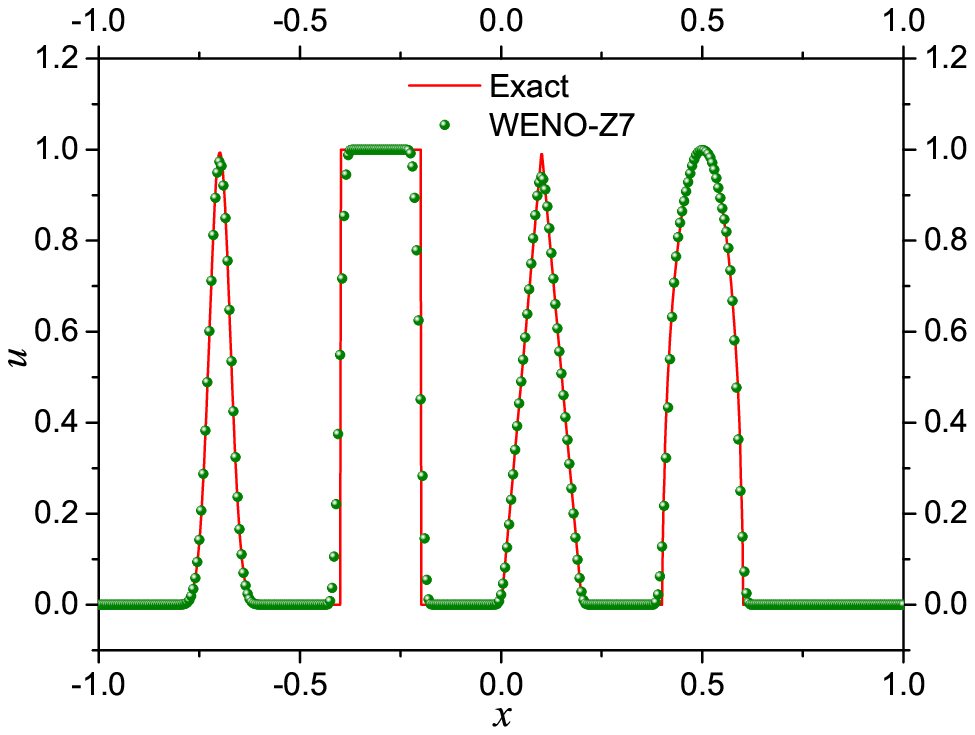}}
  \subfigure[ENO-AO7]{
  \label{FIG:Linear_Advection_ENO-AO7}
  \includegraphics[width=5.5 cm]{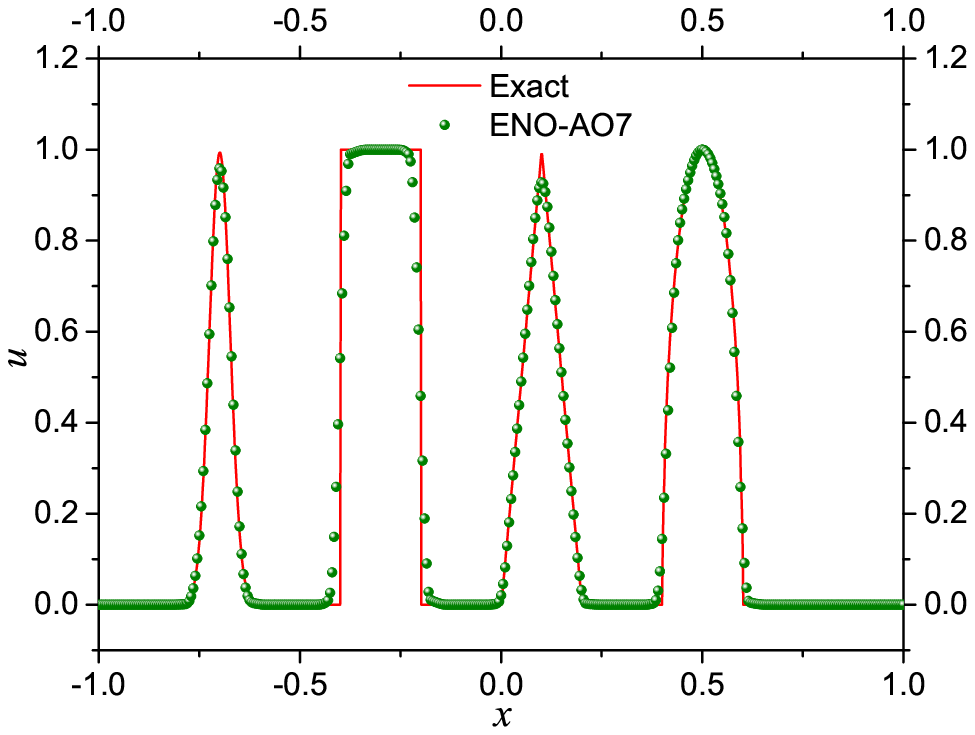}}
  \caption{The advection of a combination of Gaussians, a square wave, a sharp triangle wave,
  and a half ellipse arranged from left to right at $t=20$ calculated by WENO-Z and ENO-AO schemes with $\Delta x=1/200$.}
\label{FIG:Linear_Advection}
\end{figure}

Finally, we test the high-fidelity property of the proposed ENO-AO schemes for 
different shape of solutions by using the periodic advection of a combination of Gaussians, a square wave, a sharp
triangle wave, and a half ellipse arranged from left to right.
This case was originally proposed by Jiang and Shu \cite{Jiang_Shu1996WENO}
and was widely used to test the performance of high-order schemes.
The computational domain is $[-1,1]$,
and the initial condition is given by
\begin{subequations}
\begin{equation*}
  u_0(x)=\begin{cases}
        \frac{1}{6}\left[G(x,\beta,z-\delta)+4G(x,\beta,z)+G(x,\beta,z+\delta)\right], & \mbox{if } -0.8\le x\le -0.6 \\
        1, & \mbox{if } -0.4\le x\le -0.2 \\
        1-10\left|x-0.1\right|, & \mbox{if } 0\le x\le 0.2 \\
        \frac{1}{6}\left[F(x,\alpha,a-\delta)+4F(x,\alpha,a)+F(x,\alpha,a+\delta)\right], & \mbox{if } 0.4\le x\le 0.6 \\
        0, & \mbox{otherwise}.
  \end{cases}
\end{equation*}
\text{where,}\\
\begin{equation*}
  G(x,\beta,z)=e^{-\beta(x-z)^2},
\end{equation*}

\begin{equation*}
  F(x,\alpha,a)=\sqrt{\mbox{max}(1-\alpha^2(x-a)^2,0)}.
\end{equation*}

\end{subequations}
The constants are set as $a=0.5$, $z=-0.7$, $\delta=0.005$, $\alpha=10$, and $\beta=\frac{\log2}{36\delta^2}$.
Periodic boundary conditions are implemented on the left and right sides.
Fig. \ref{FIG:Linear_Advection} shows the profile of $u$ 
at $t=20$ calculated by WENO-Z and ENO-AO schemes with $\Delta x=1/200$.
We observe that ENO-AO schemes are as accuracy as WENO-Z schemes for all shape of solutions,
although they contains lower-order stencils than WENO-Z schemes.
In fact, if we oberve carefully, ENO-AO5 performs better than WENO-Z5 for the smooth Gaussians.
This case demonstrates that the new smoothness indicators can correctly measure
the smoothness of stencils with unequal sizes and can help to select the optimal stencil
for different solutions.

\subsection{Numerical examples for the 1D Euler equations}
We consider the 1D compressible Euler equations
\begin{equation}
  \frac{\partial \mathbf{U}}{\partial t}+\frac{\partial \mathbf{F}}{\partial x}=0,
\end{equation}
where $\mathbf{U}=[\rho,\rho u,\rho e]^T$ and $\mathbf{F}=[{\rho u,\rho u^2+p,(\rho e+p)u}]^T$
with $\rho$, $u$, $p$, and $e$ denoting the density, velocity, pressure, and specific total energy respectively.
We assume the gas is ideal and the specific total energy is calculated as $e=\frac{p}{\rho(\gamma-1)}+\frac{1}{2} u^2$.
where $\gamma$ is the specific heat ratio.

\begin{figure}
  \centering
  \subfigure[]{
  \label{FIG:Lax_Density_WENO_Z5}
  \includegraphics[width=5.5 cm]{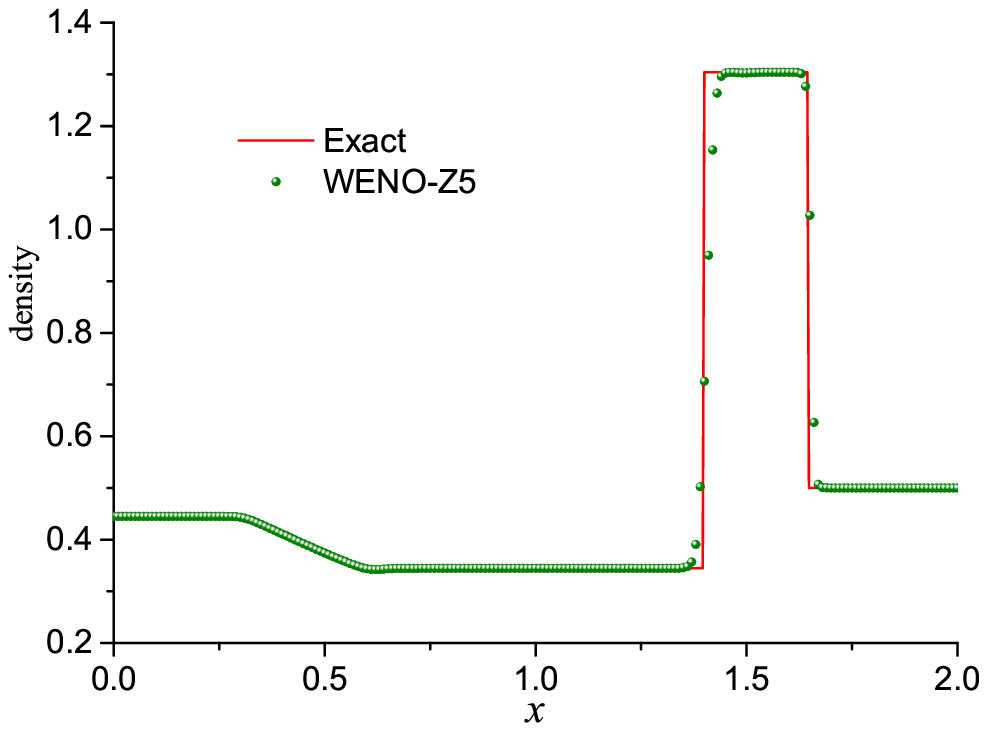}}
  \subfigure[]{
  \label{FIG:Lax_Velocity_WENO_Z5}
  \includegraphics[width=5.5 cm]{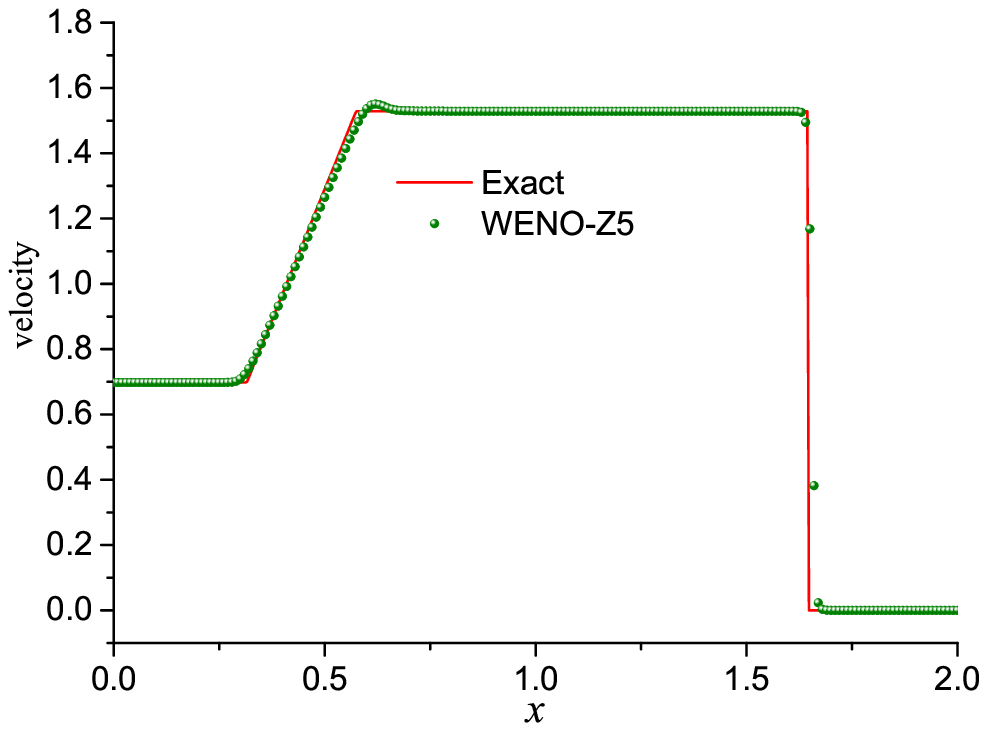}}
  \subfigure[]{
  \label{FIG:Lax_Density_ENO_AO5}
  \includegraphics[width=5.5 cm]{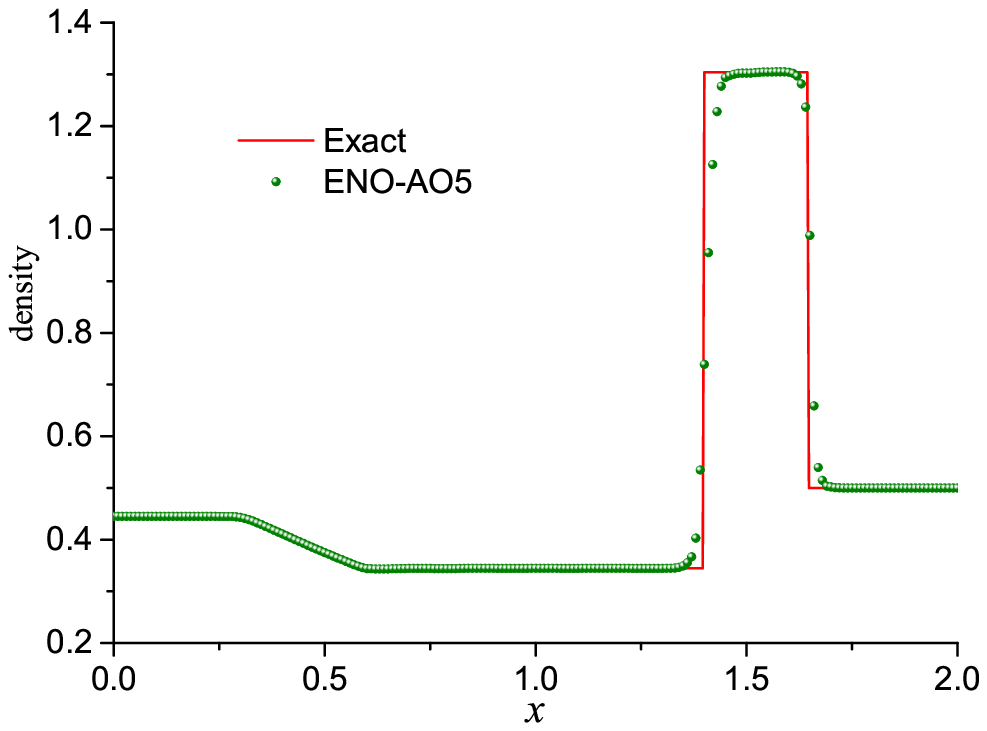}}
  \subfigure[]{
  \label{FIG:Lax_Velocity_ENO_AO5}
  \includegraphics[width=5.5 cm]{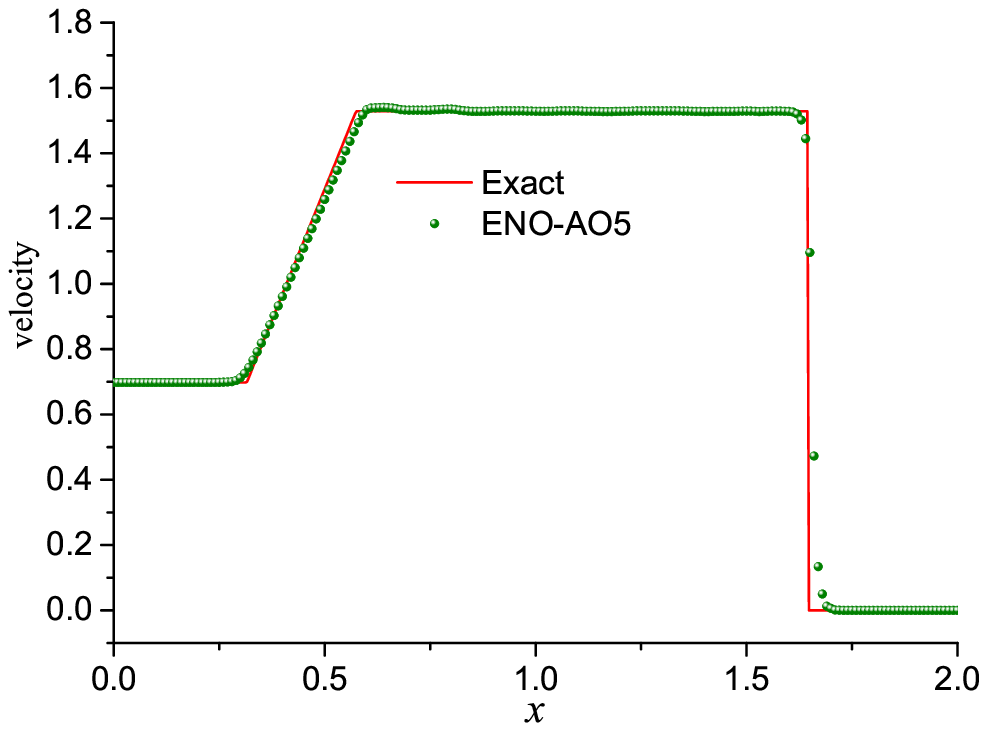}}
  \subfigure[]{
  \label{FIG:Lax_Density_WENO_Z7}
  \includegraphics[width=5.5 cm]{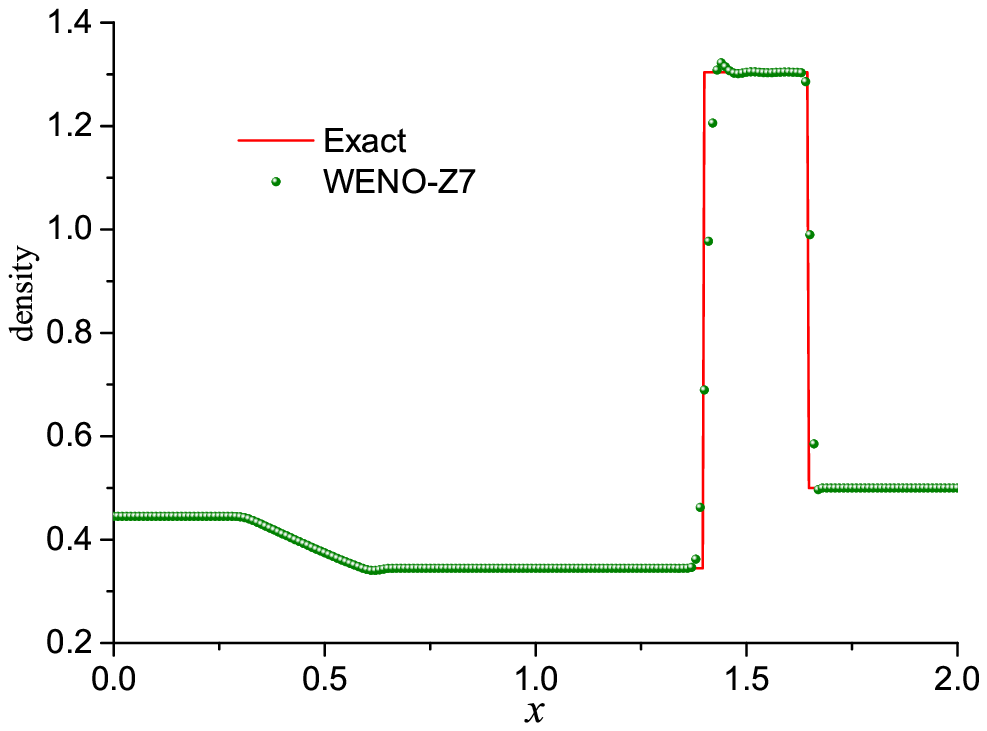}}
  \subfigure[]{
  \label{FIG:Lax_Velocity_WENO_Z7}
  \includegraphics[width=5.5 cm]{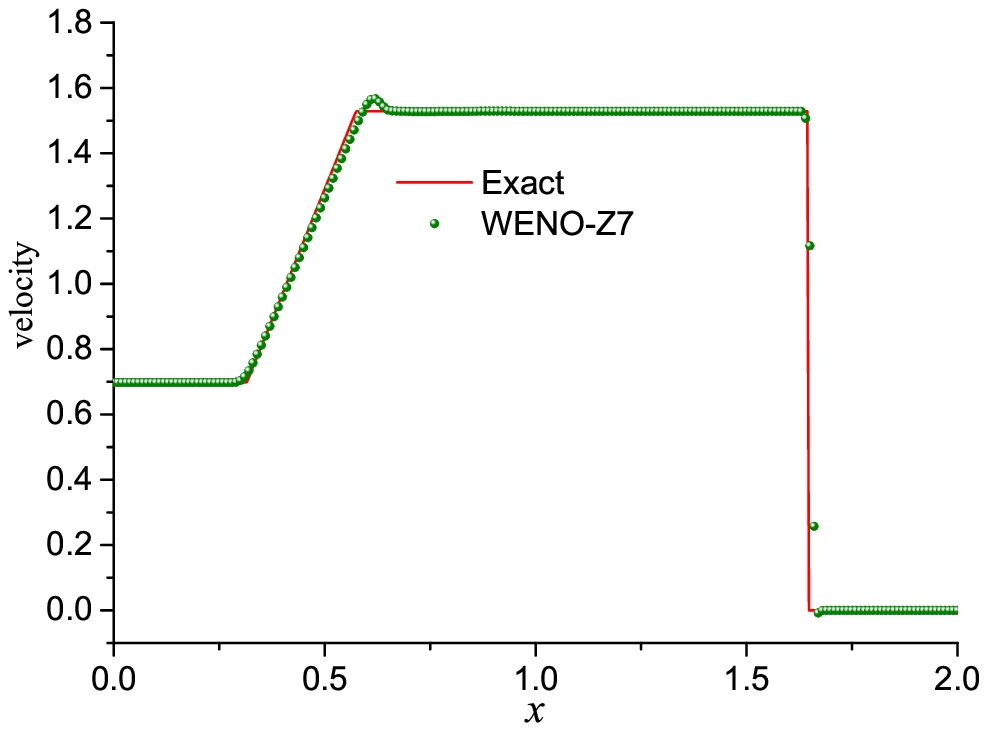}}
  \subfigure[]{
  \label{FIG:Lax_Density_ENO_AO7}
  \includegraphics[width=5.5 cm]{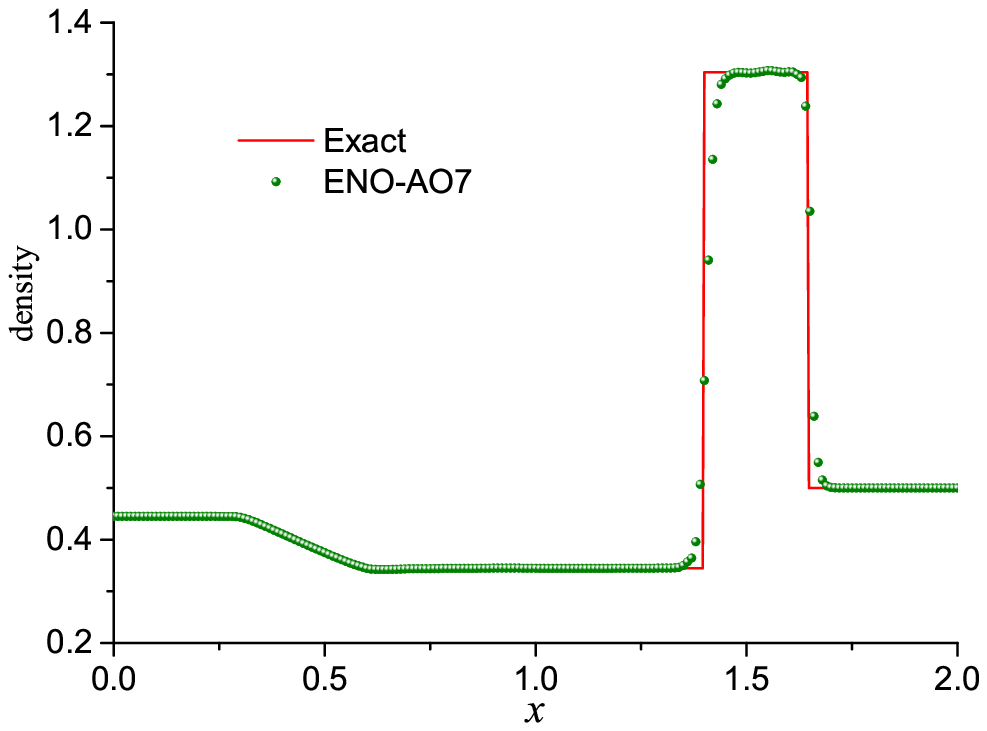}}
  \subfigure[]{
  \label{FIG:Lax_Velocity_ENO_AO7}
  \includegraphics[width=5.5 cm]{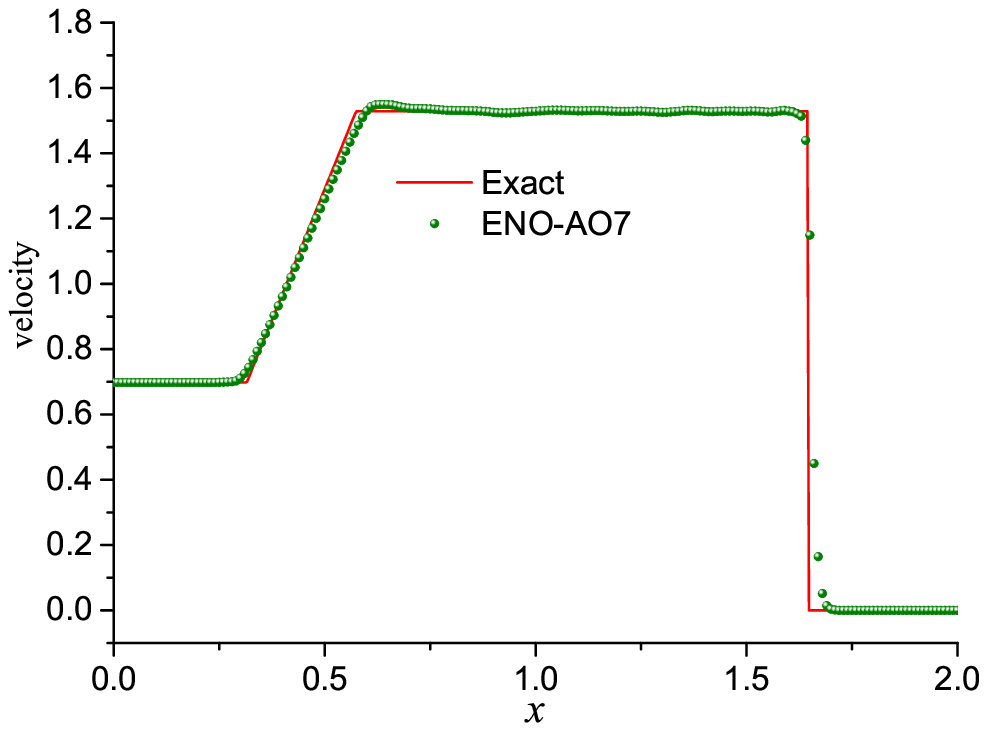}}
  \caption{The density profiles of the Lax shock tube problem at $t=0.26$ calculated by WENO-Z and ENO-AO with $\Delta x=1/100$.}
\label{FIG:Lax}
\end{figure}

The first case is the Lax \cite{Lax1954} shock tube problem that is used to test the shock-capturing property of numerical methods.
The computational domain is [0,2], the ratio of specific heats $\gamma=1.4$, and the initial condition is given by
\begin{equation*}
  (\rho,u,p)=\begin{cases}
               (0.445,0.698,3.528), & \mbox{if } x<1 \\
               (0.5,0,0.571), & \mbox{otherwise}.
             \end{cases}
\end{equation*}
Non-reflection boundary conditions are applied on the left and right boundaries.
Fig. \ref{FIG:Lax} shows the density and velocity profiles at $t=0.26$ 
calculated by WENO-Z and ENO-AO with $\Delta x=1/100$.
We observe that the velocity profiles of WENO-Z schemes have a spike near the expansion wave,
and the density profile of the WENO-Z7 scheme is also oscillatory near the contact discontinuity.
The corresponding ENO-AO schemes are more robust since they contain low-order stencils 
which are suitable to capture discontinuities.

\begin{figure}
  \centering
  \subfigure[]{
  \label{FIG:Titarev_Density_WENO_Z5}
  \includegraphics[width=5.5 cm]{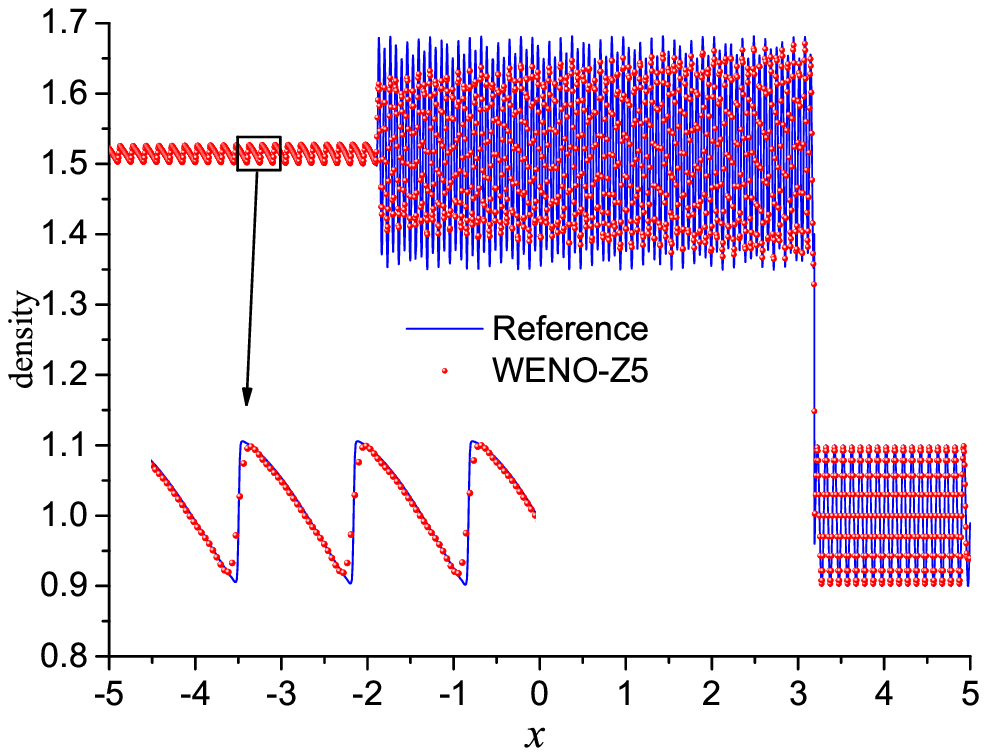}}
  \subfigure[]{
  \label{FIG:Titarev_Density_WENO_Z5_enlarge}
  \includegraphics[width=5.5 cm]{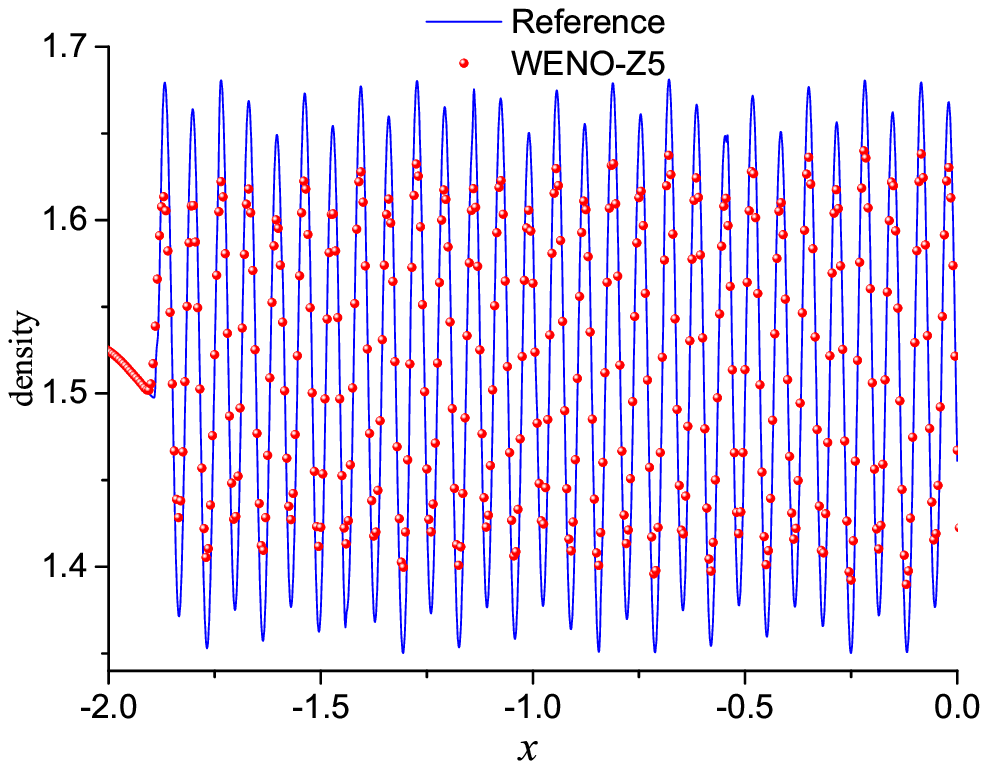}}
  \subfigure[]{
  \label{FIG:Titarev_Density_ENO_AO5}
  \includegraphics[width=5.5 cm]{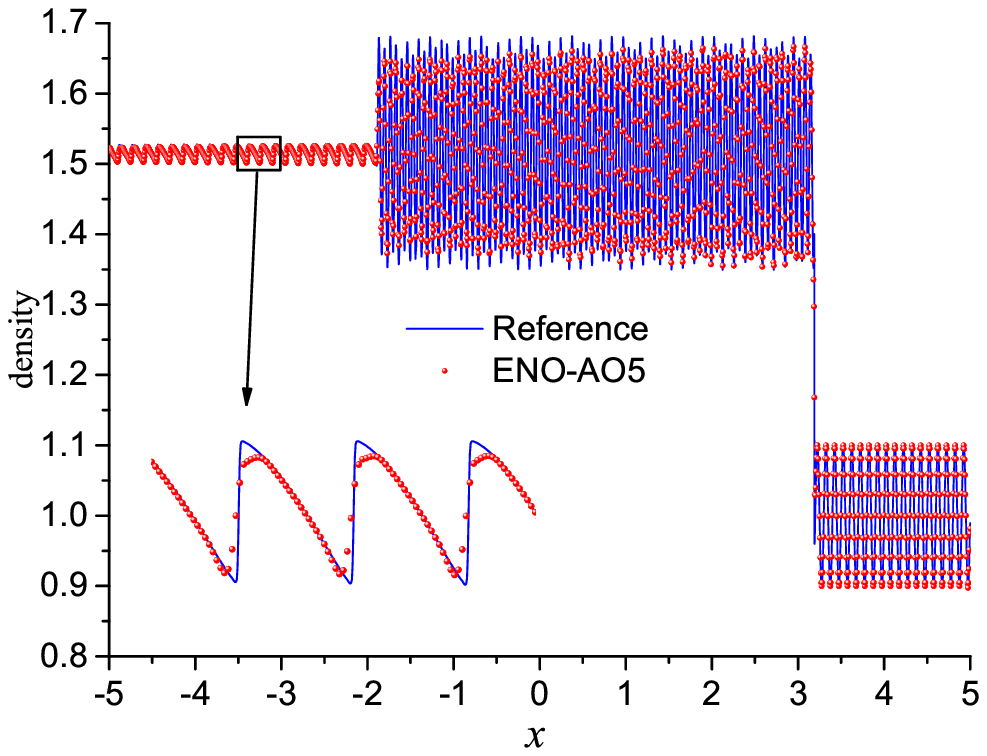}}
  \subfigure[]{
  \label{FIG:Titarev_Density_ENO_AO5_enlarge}
  \includegraphics[width=5.5 cm]{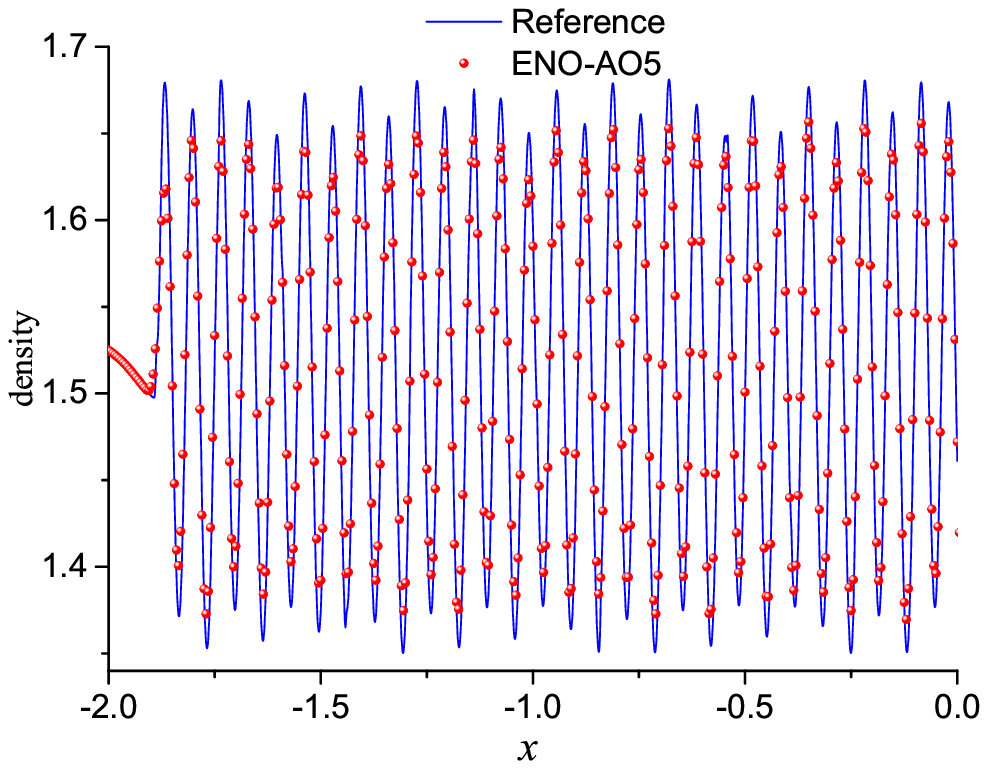}}
  \subfigure[]{
  \label{FIG:Titarev_Density_WENO_Z7}
  \includegraphics[width=5.5 cm]{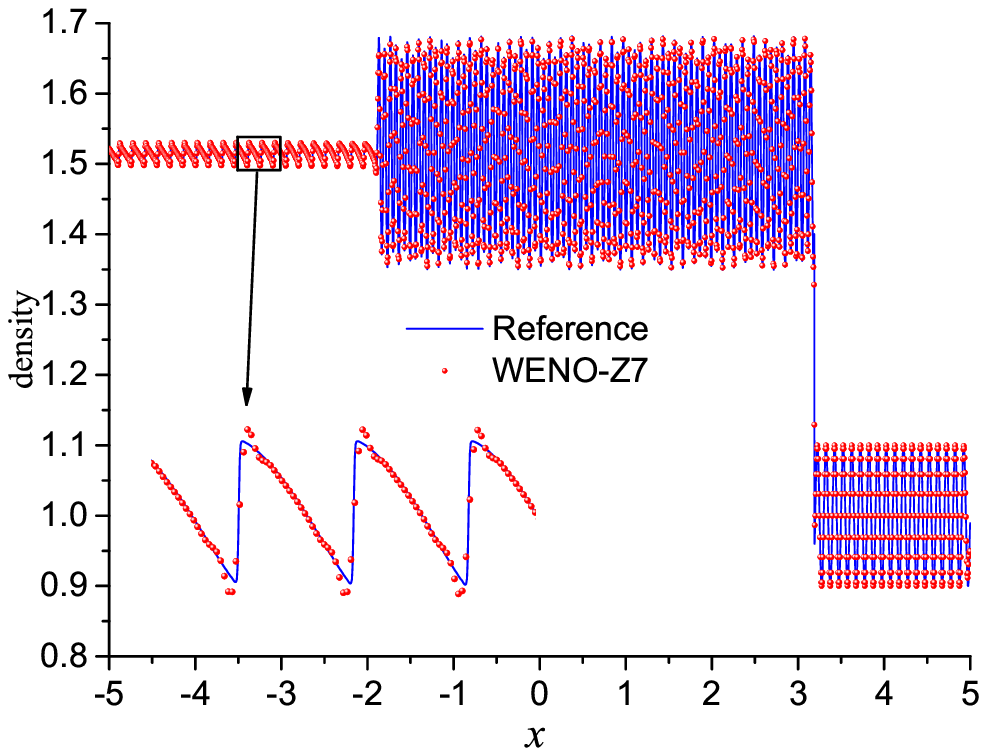}}
  \subfigure[]{
  \label{FIG:Titarev_Density_WENO_Z7_enlarge}
  \includegraphics[width=5.5 cm]{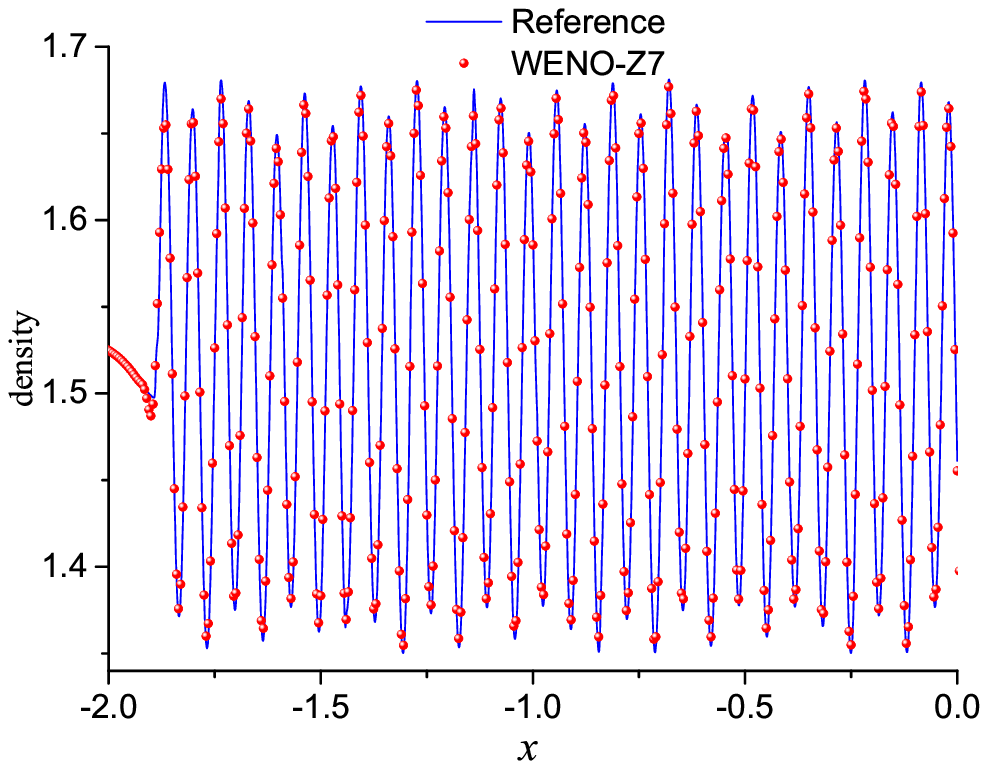}}
  \subfigure[]{
  \label{FIG:Titarev_Density_ENO_AO7}
  \includegraphics[width=5.5 cm]{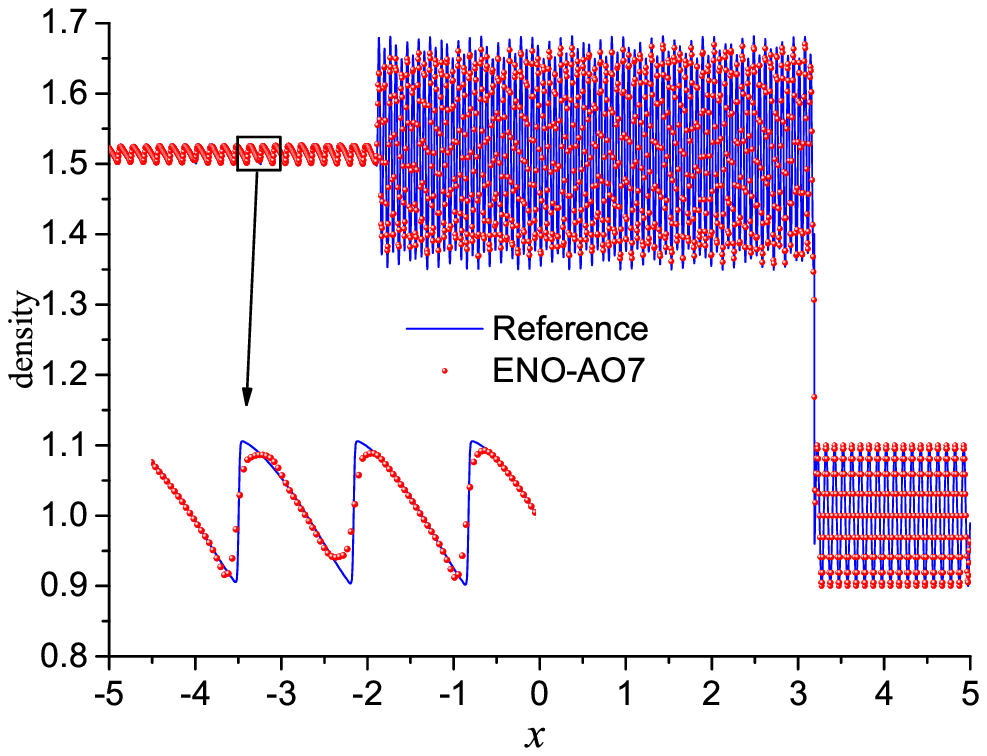}}
  \subfigure[]{
  \label{FIG:Titarev_Density_ENO_AO7_enlarge}
  \includegraphics[width=5.5 cm]{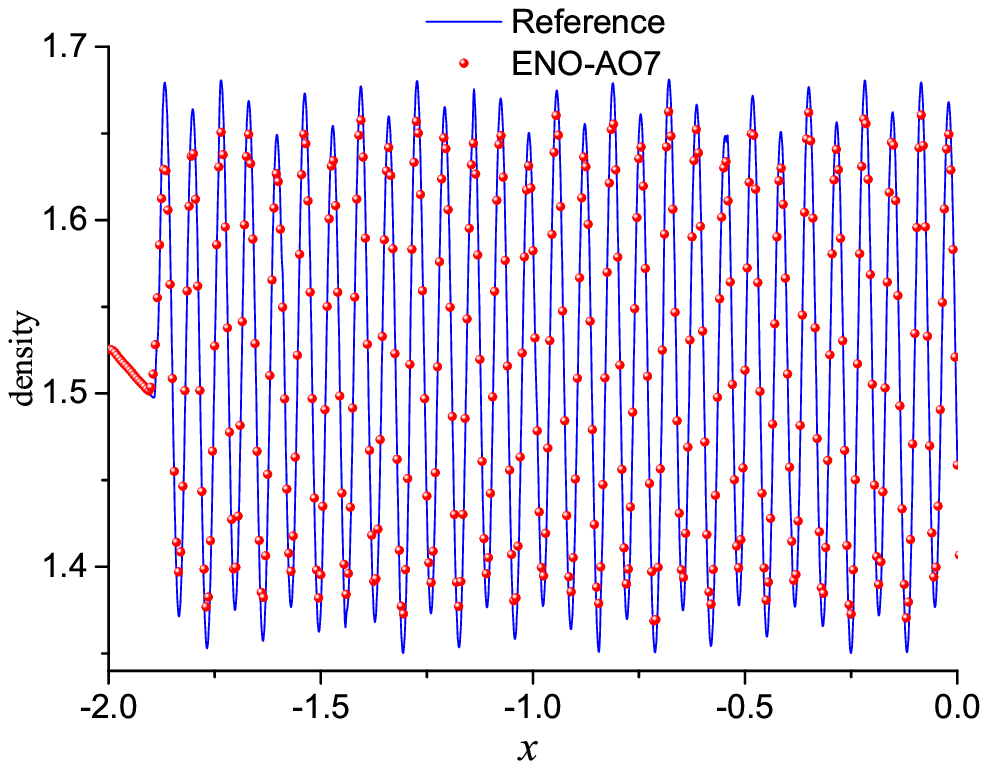}}
  \caption{The density profiles of the Titarev and Toro problem at $t=5$ calculated by WENO-Z and ENO-AO with $\Delta x=1/200$.}
\label{FIG:Titarev_Toro}
\end{figure}

The second case is the Titarev and Toro\cite{Titarev2004WENO_FV} problem 
which is an upgraded version of the Shu-Osher problem \cite{Shu1988EfficientENO}.
It depicts the interaction of a Mach 1.1 moving shock with a high-frequency entropy sine wave 
which contains very complicated flow structures.
The computational domain is [-5,5], the ratio of specific heats $\gamma=1.4$, and the initial condition is given by
\begin{equation*}
  (\rho,u,p)=\begin{cases}
               (1.515695,0.523346,1.805), & \mbox{if } x<-4.5 \\
               (1+0.1\mbox{sin}(20\pi x),0,1), & \mbox{otherwise}.
             \end{cases}
\end{equation*}
Non-reflection boundary conditions are applied on the left and right sides.
Fig. \ref{FIG:Titarev_Toro} shows the density profiles at $t=5$ 
calculated by WENO-Z and ENO-AO schemes with $\Delta x=1/200$.
The reference solution is calculated by the WENO-Z5 scheme with $\Delta x=1/1000$.
We observe that the result of ENO-AO5 is closer to the reference solution in the high-frequency region than that of WENO-Z5.
The result of WENO-Z7 is closest to the reference solution in $[-2,3]$, 
but it has some overshoots and undershoots in $[0,-2]$, as shown in Fig. \ref{FIG:Titarev_Density_WENO_Z7}.

\subsection{Numerical examples for the 2D Euler equations}
The 2D compressible Euler equations can be written as 
\begin{equation}\label{Eq:2D_Euler_equations}
  \frac{\partial \mathbf{U}}{\partial t}+\frac{\partial \mathbf{F}}{\partial x}+\frac{\partial \mathbf{G}}{\partial y}=0,
\end{equation}
with $\mathbf{U}=[\rho,\rho u,\rho v,\rho e]^T$, $\mathbf{F}=[{\rho u,\rho u^2+p,\rho uv,(\rho e+p)u}]^T$,
and $\mathbf{G}=[{\rho v,\rho uv,\rho v^2+p,(\rho e+p)v}]^T$,
where $\rho$, $u$, $v$, $p$, and $e$ denote the density, $x$-velocity, $y$-velocity, pressure, and specific total energy respectively.
We use the ideal gas model and the specific total energy is calculated as $e=\frac{p}{\rho(\gamma-1)}+\frac{1}{2} (u^2+v^2)$.

\subsubsection{Two-dimensional Riemann problems}
\begin{figure}
  \centering
  \subfigure[WENO-Z5]{
  \label{FIG:RP1_Density_WENO_Z5}
  \includegraphics[width=5.5 cm]{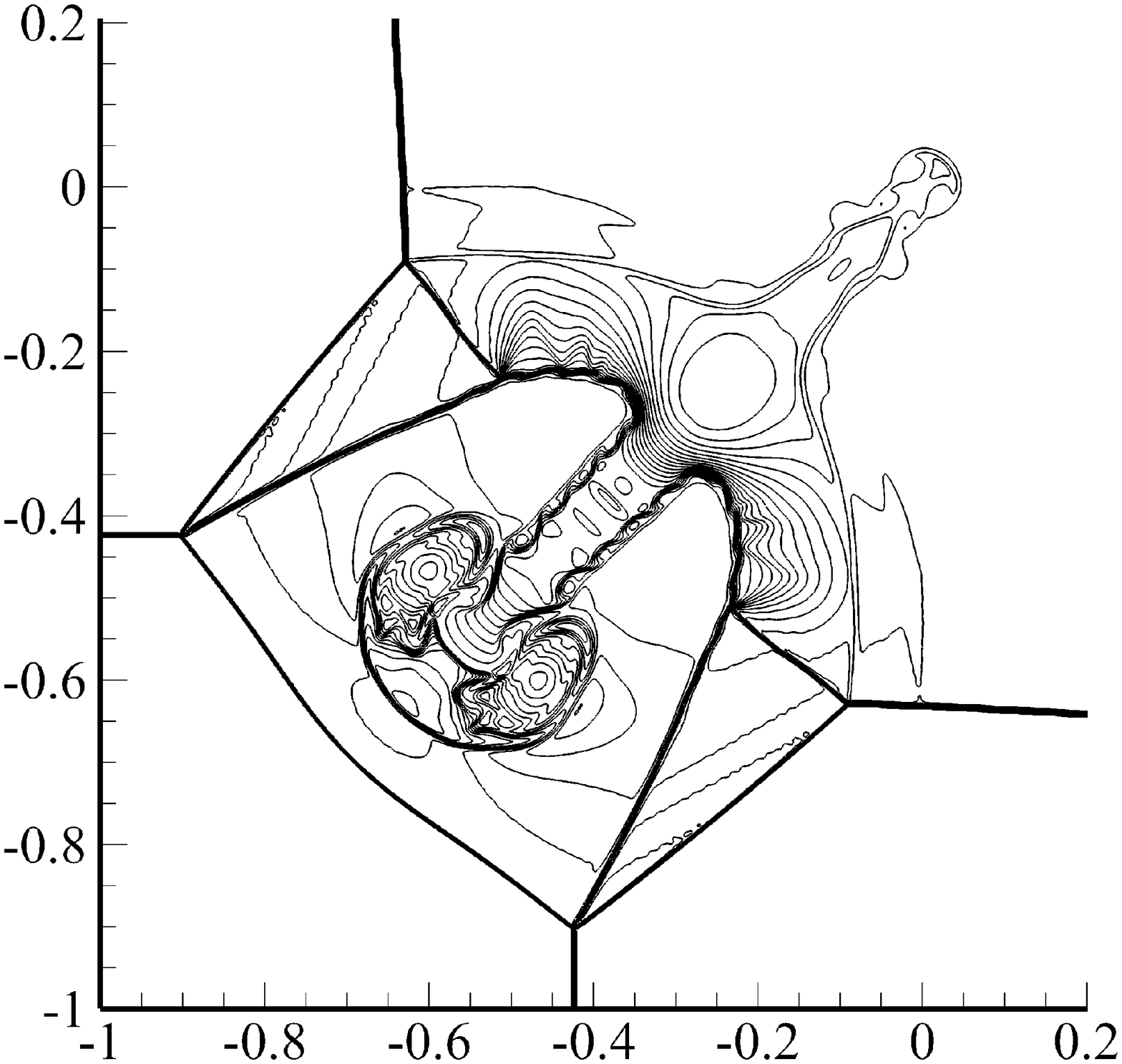}}
  \subfigure[ENO-AO5]{
  \label{FIG:RP1_Density_ENO_AO5}
  \includegraphics[width=5.5 cm]{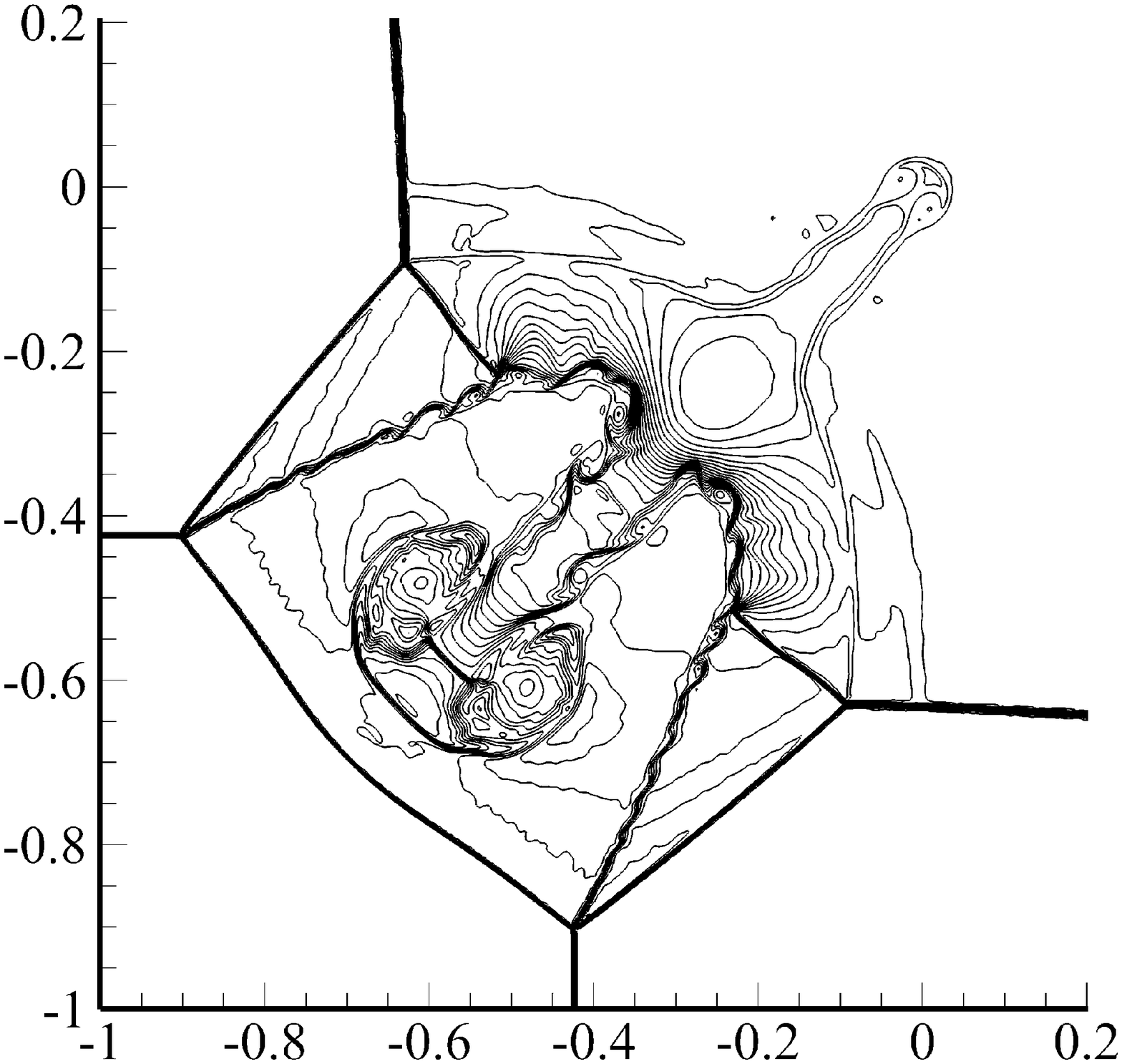}}
  \subfigure[WENO-Z7]{
  \label{FIG:RP1_Density_WENO_Z7}
  \includegraphics[width=5.5 cm]{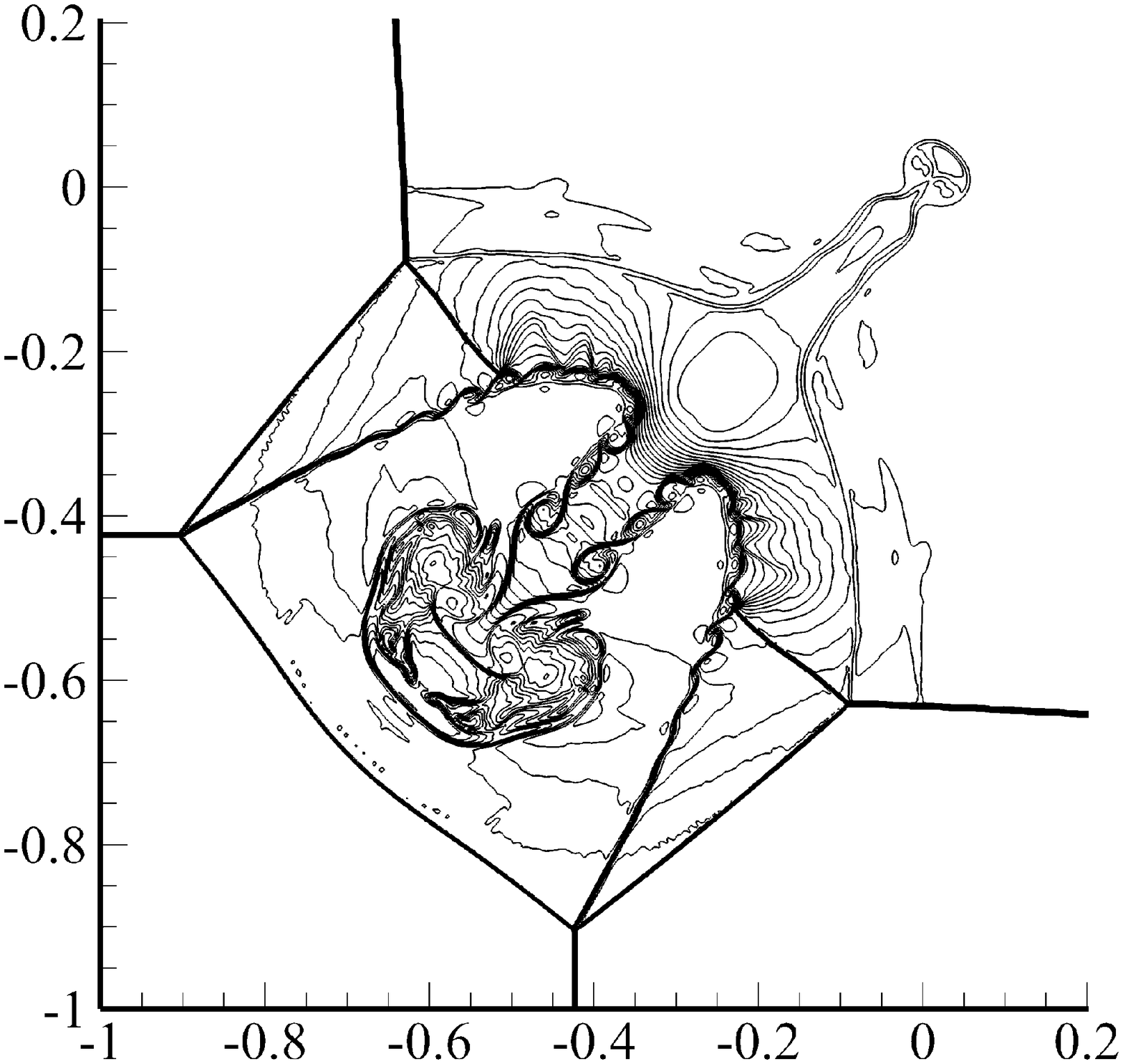}}
  \subfigure[ENO-AO7]{
  \label{FIG:RP1_Density_ENO_AO7}
  \includegraphics[width=5.5 cm]{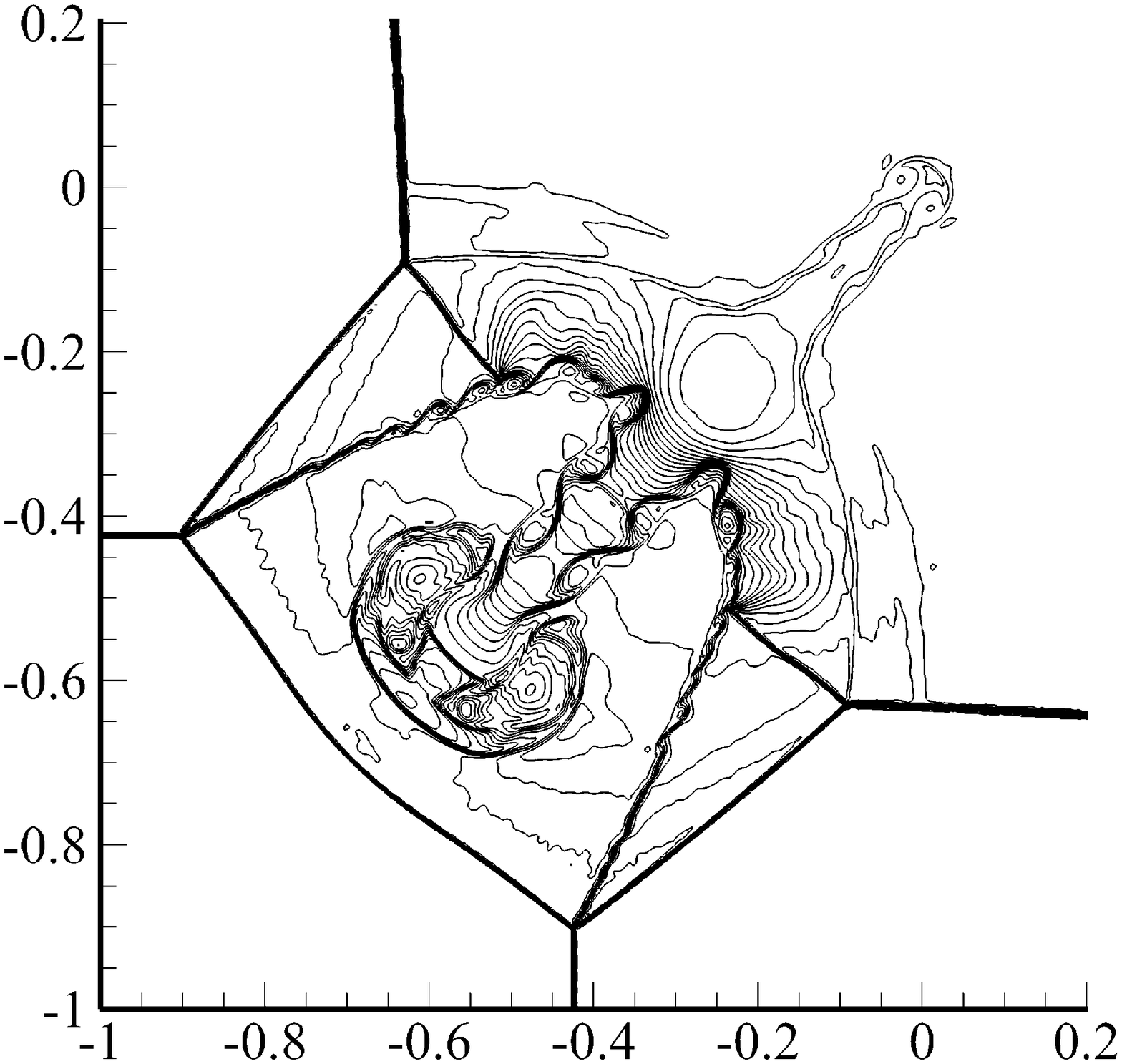}}
  \caption{The density contours of the first 2D Riemann problem at $t=1$ calculated 
  by WENO-Z and ENO-AO with $\Delta x=\Delta y=1/400$.
  The density contours contain 30 equidistant contours from 0.2 to 1.7.}
\label{FIG:RP1}
\end{figure}

We consider two different configurations with the ratio of specific heats $\gamma=1.4$.
The initial condition of the first configuration is given by
\begin{equation*}
  \left(\rho, u, v, p\right)=
  \begin{cases}
    \left(0.138,1.206,1.206,0.029\right),    & \text{in }[-1, 0]\times[-1, 0],\\
    \left(0.5323,1.206,0,0.3\right),         & \text{in }[-1, 0]\times[0, 1],\\
    \left(1.5,0,0,1.5\right),                & \text{in }[0, 1] \times[0, 1],\\
    \left(0.5323,0,1.206,0.3\right),         & \text{in }[0, 1] \times[-1, 0],
  \end{cases}
\end{equation*}
which depicts two horizontally moving shocks and two vertically moving shocks.
When the simulation begins, the interaction of the four normal shocks will 
result in two double-Mach reflections and an oblique shock moving along the diagonal of the computational domain.
The density contours at $t=1$ calculated by WENO-Z and ENO-AO schemes with $801\times801$ mesh points 
are shown in Fig. \ref{FIG:RP1}.
In all simulations, Kelvin–Helmholtz (KH) instabilities occur along the slip lines.
ENO-AO5 captures more details than WENO-Z5,
and the results computed by ENO-AO7 and WENO-Z7 are comparable.

\begin{figure}
  \centering
  \subfigure[WENO-Z5]{
  \label{FIG:RP2_Density_WENO_Z5}
  \includegraphics[width=5.5 cm]{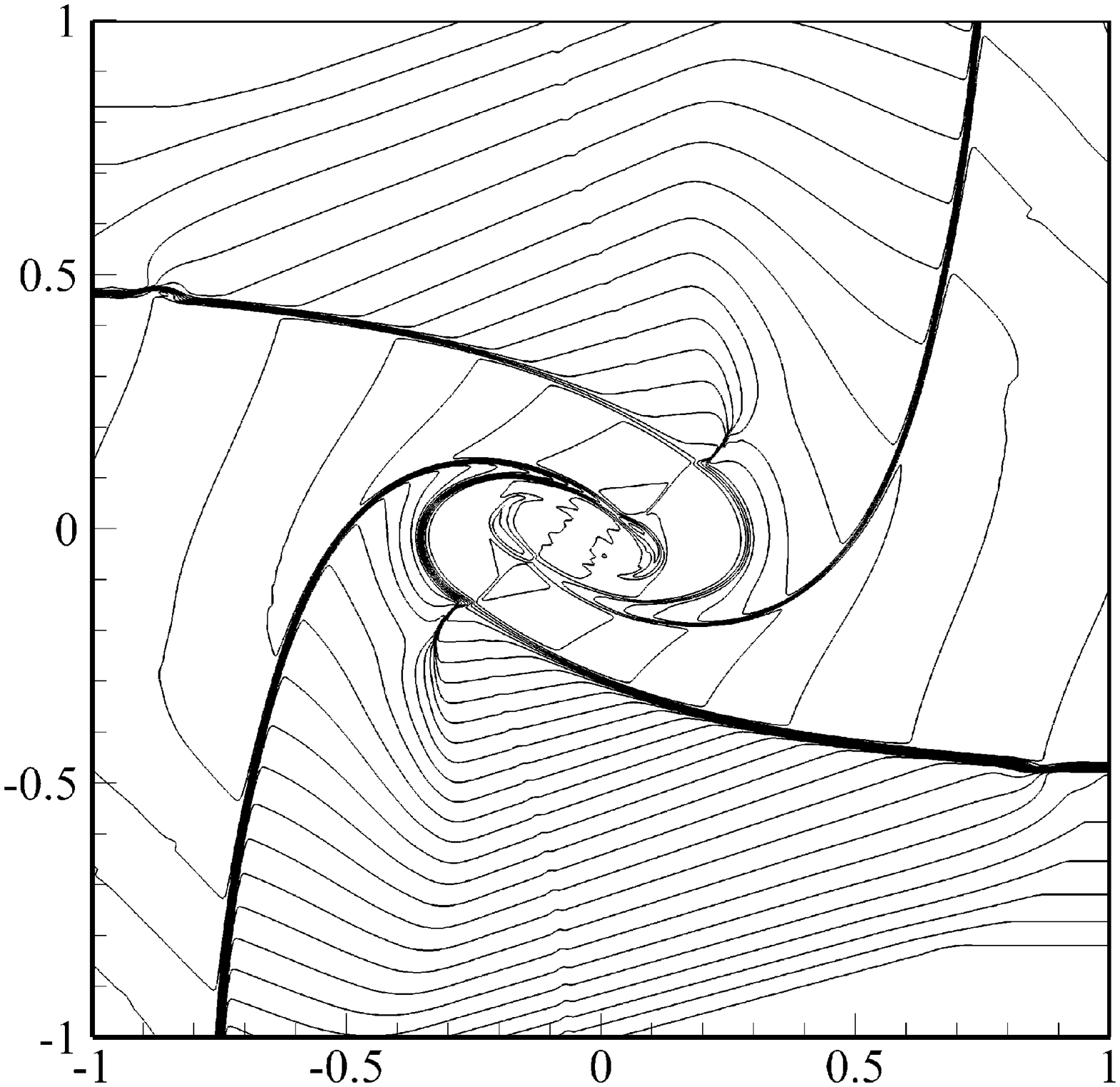}}
  \subfigure[ENO-AO5]{
  \label{FIG:RP2_Density_ENO_AO5}
  \includegraphics[width=5.5 cm]{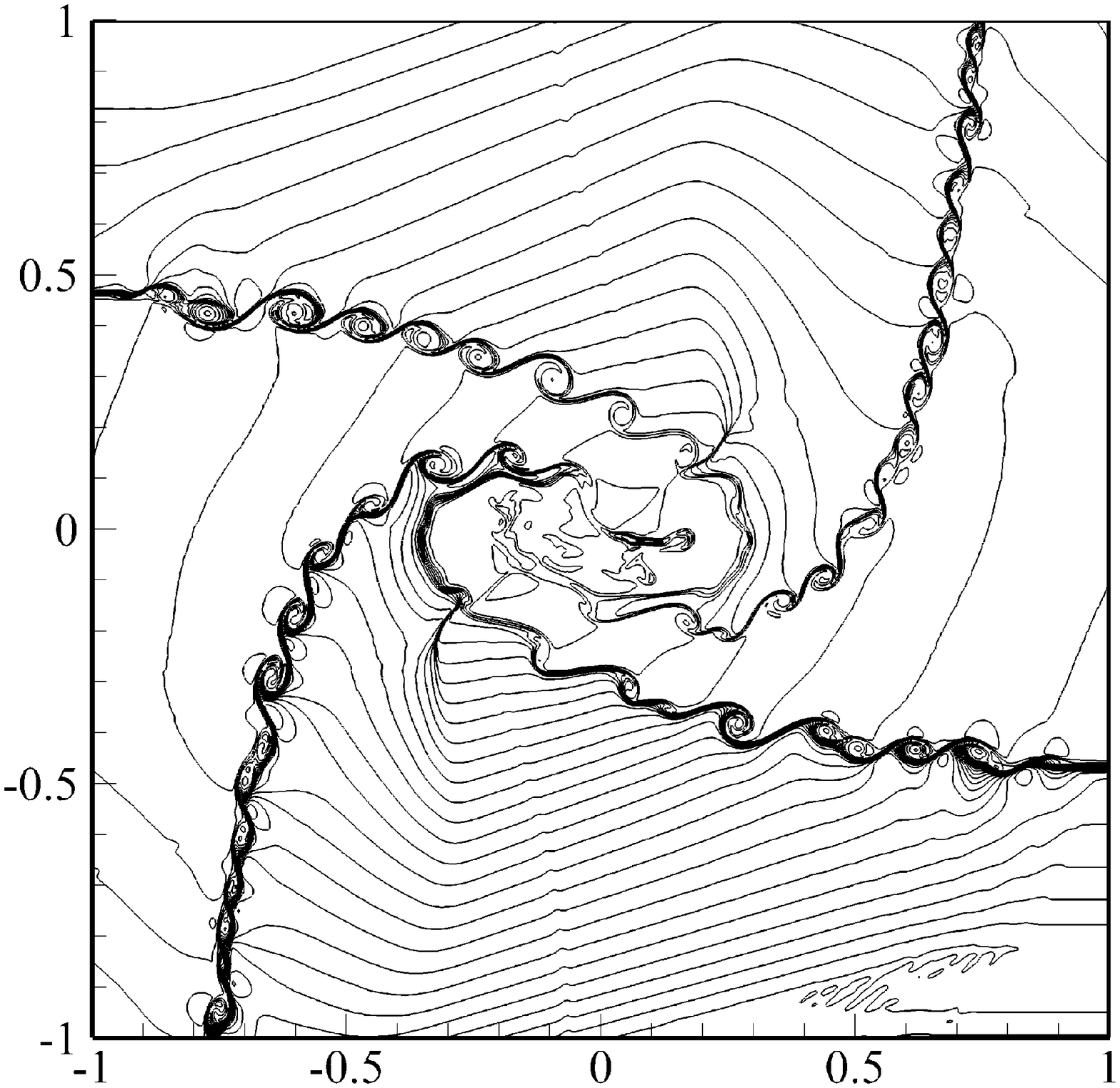}}
  \subfigure[WENO-Z7]{
  \label{FIG:RP2_Density_WENO_Z7}
  \includegraphics[width=5.5 cm]{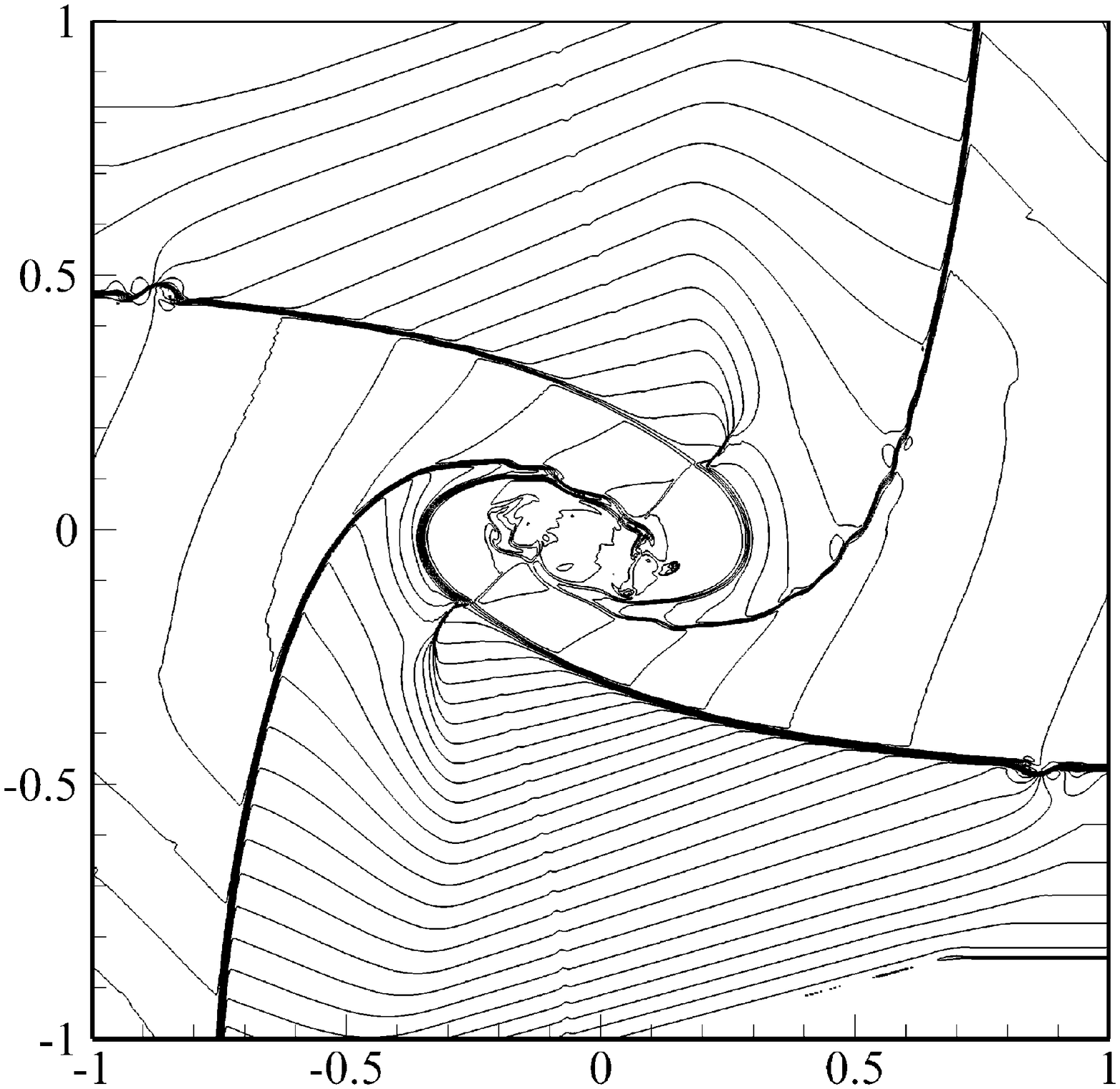}}
  \subfigure[ENO-AO7]{
  \label{FIG:RP2_Density_ENO_AO7}
  \includegraphics[width=5.5 cm]{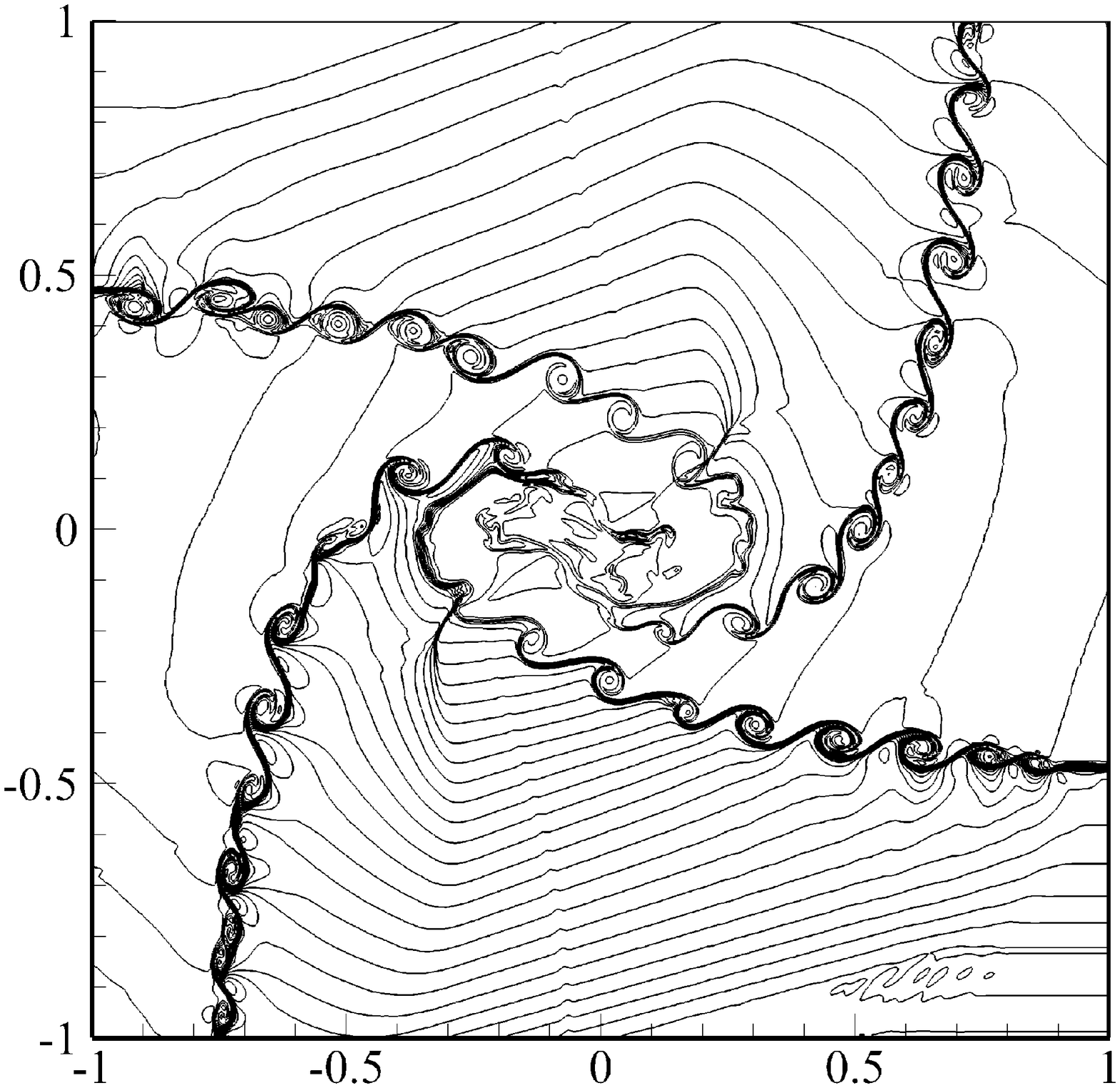}}
  \caption{The density contours of the second 2D Riemann problem at $t=1$ calculated 
  by WENO-Z and ENO-AO with $\Delta x=\Delta y=1/400$.
  The density contours contain 30 equidistant contours from 0.2 to 3.}
\label{FIG:RP2}
\end{figure}

The second configuration starts off four contact discontinuities,
and the initial condition is given by
\begin{equation*}
  \left(\rho, u, v, p\right)=
  \begin{cases}
    \left(1,-0.75,0.5,1\right),        & \text{in }[-1, 0]\times[-1, 0],\\
    \left(2,0.75,0.5,1\right),         & \text{in }[-1, 0]\times[0, 1],\\
    \left(1,0.75,-0.5,1\right),        & \text{in }[0, 1] \times[0, 1],\\
    \left(3,-0.75,-0.5,1\right),       & \text{in }[0, 1] \times[-1, 0].
  \end{cases}
\end{equation*}
Fig. \ref{FIG:RP2} shows the density contours at $t=1$ 
calculated by WENO-Z and ENO-AO schemes with $801\times801$ mesh points.
All the schemes capture the large-scale spiral at the center of the domain,
but ENO-AO schemes capture the small-scale vortices along contact discontinuities
which are absent in the results of WENO-Z schemes.
For this configuration, ENO-AO schemes have a prominent advantage over WENO-Z schemes.

\subsubsection{Double Mach reflection problem}
\begin{figure}
  \centering
  \subfigure[WENO-Z5]{
  \label{FIG:DMR_Density_WENO_Z5}
  \includegraphics[width=10 cm]{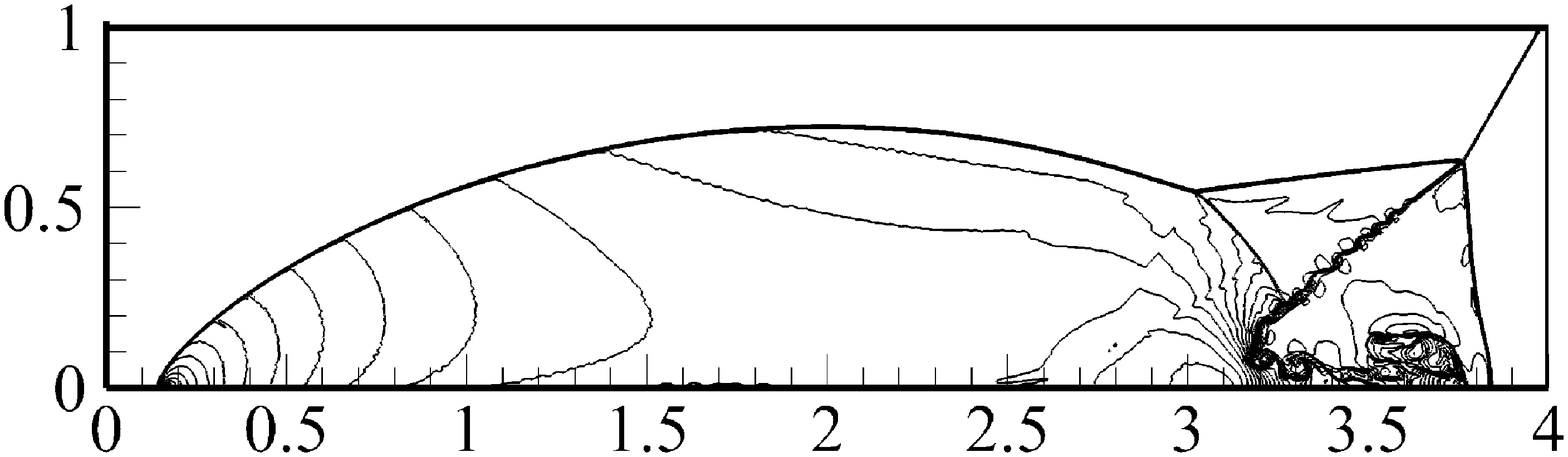}}
  \subfigure[ENO-AO5]{
  \label{FIG:DMR_Density_ENO_AO5}
  \includegraphics[width=10 cm]{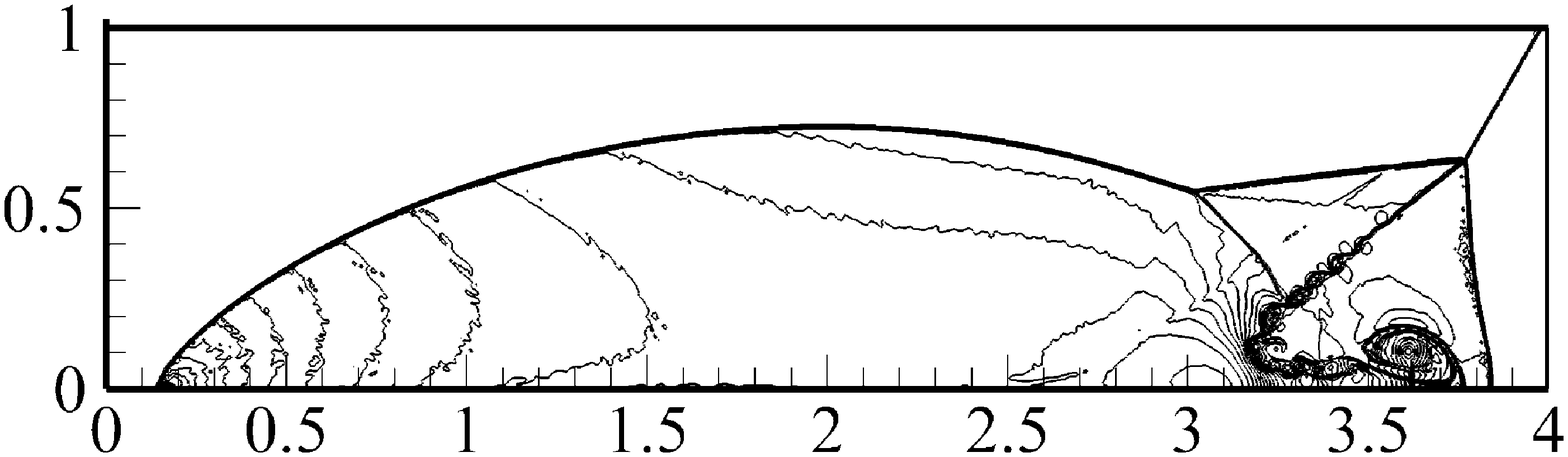}}
  \subfigure[WENO-Z7]{
  \label{FIG:DMR_Density_WENO_Z7}
  \includegraphics[width=10 cm]{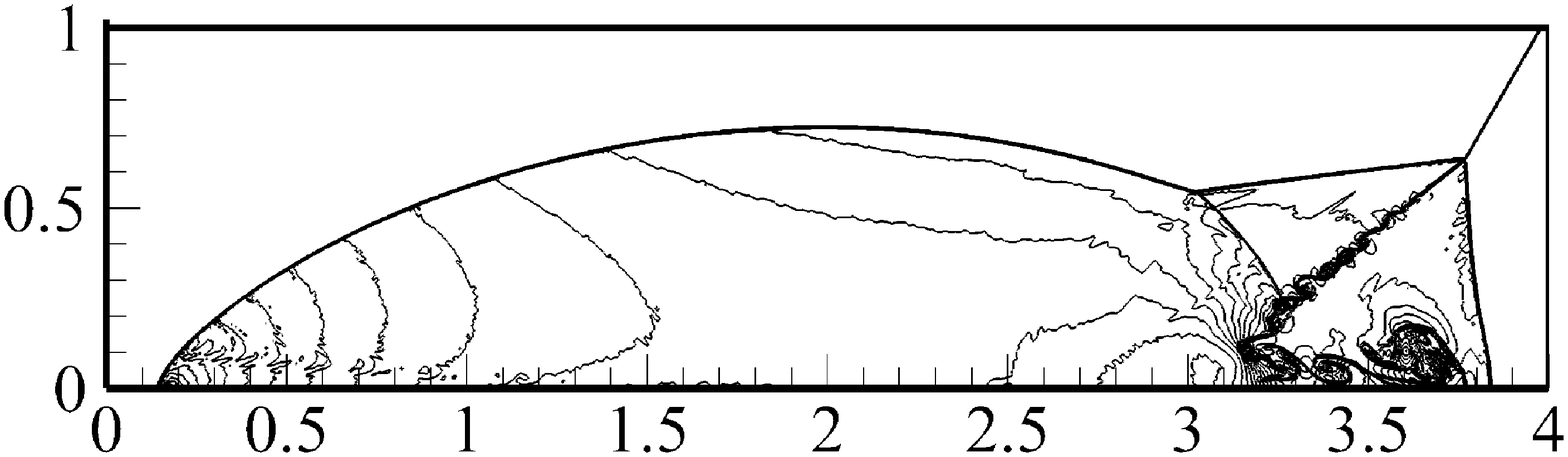}}
  \subfigure[ENO-AO7]{
  \label{FIG:DMR_Density_ENO_AO7}
  \includegraphics[width=10 cm]{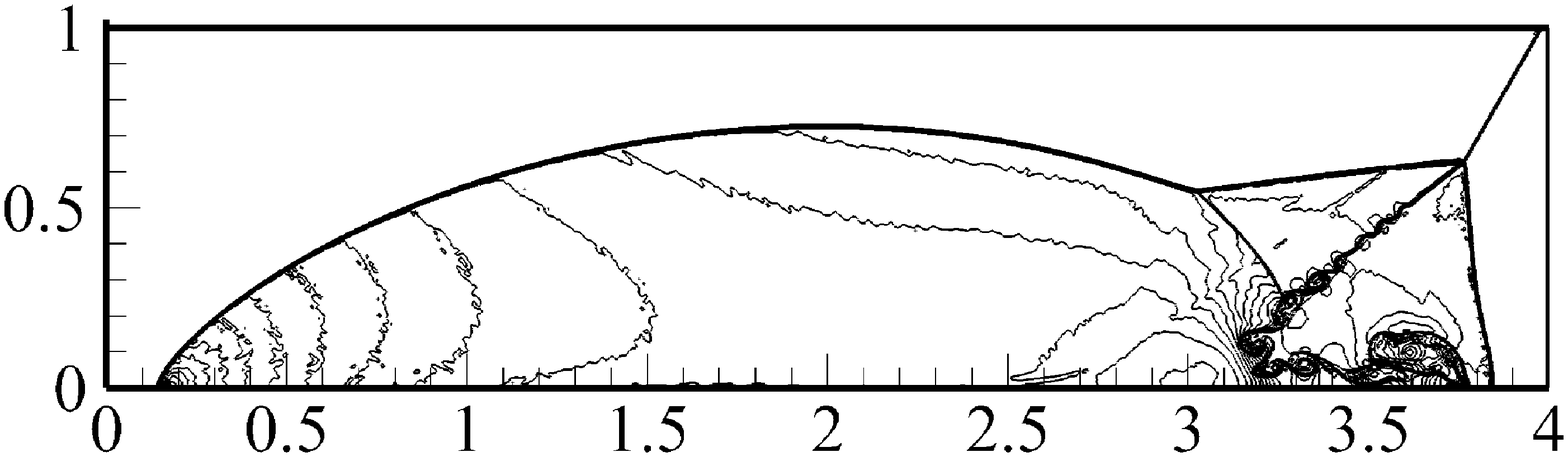}}
  \caption{The density contours of the double Mach reflection problem at $t=0.28$ calculated 
  by WENO-Z and ENO-AO with $\Delta x=\Delta y=1/400$.
  The density contours contain 40 equidistant contours from 2 to 22.}
\label{FIG:DMR}
\end{figure}

\begin{figure}
  \centering
  \subfigure[WENO-Z5]{
  \label{FIG:DMR_Density_WENO_Z5_enlarge}
  \includegraphics[width=5.5 cm]{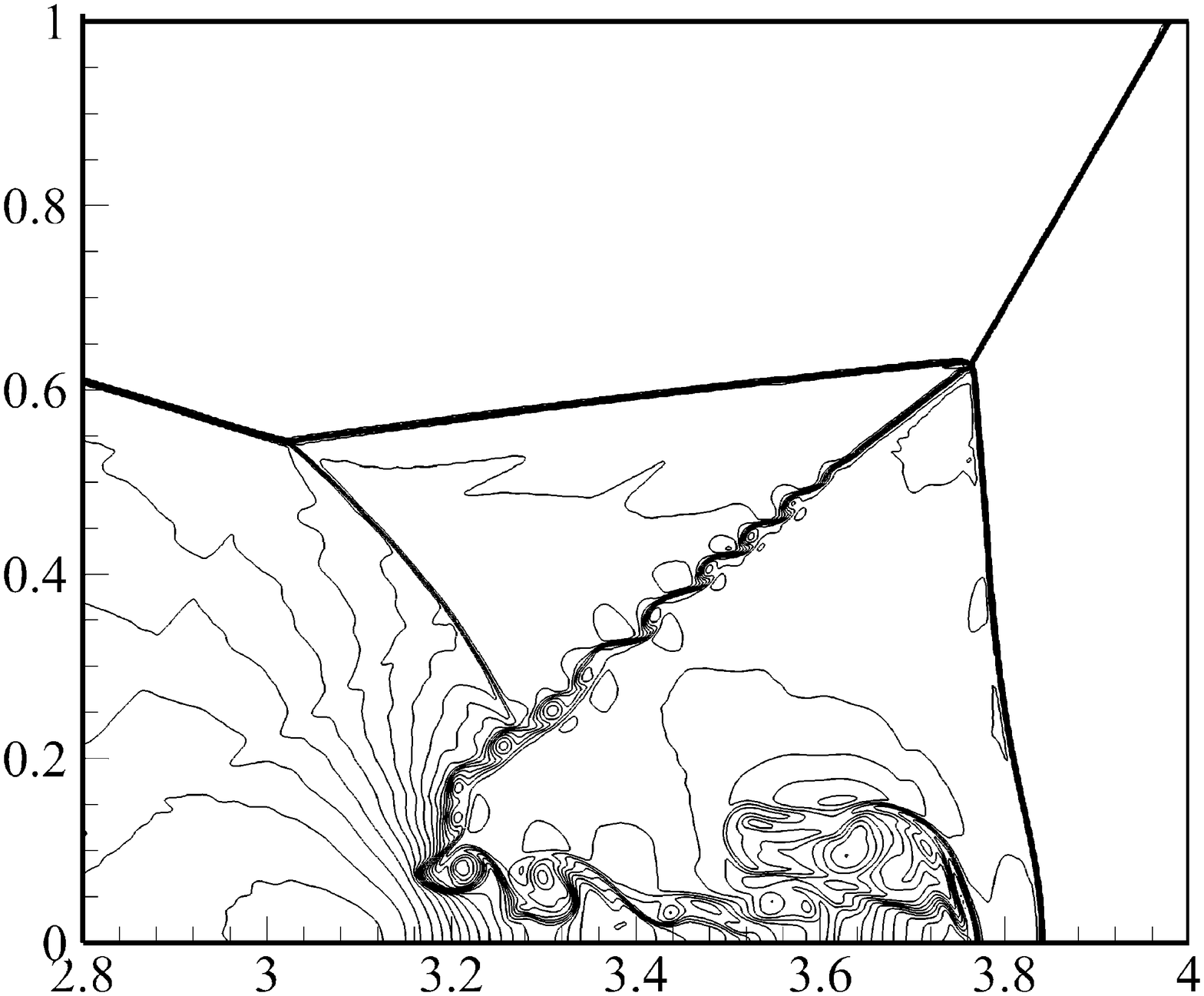}}
  \subfigure[ENO-AO5]{
  \label{FIG:DMR_Density_ENO_AO5_enlarge}
  \includegraphics[width=5.5 cm]{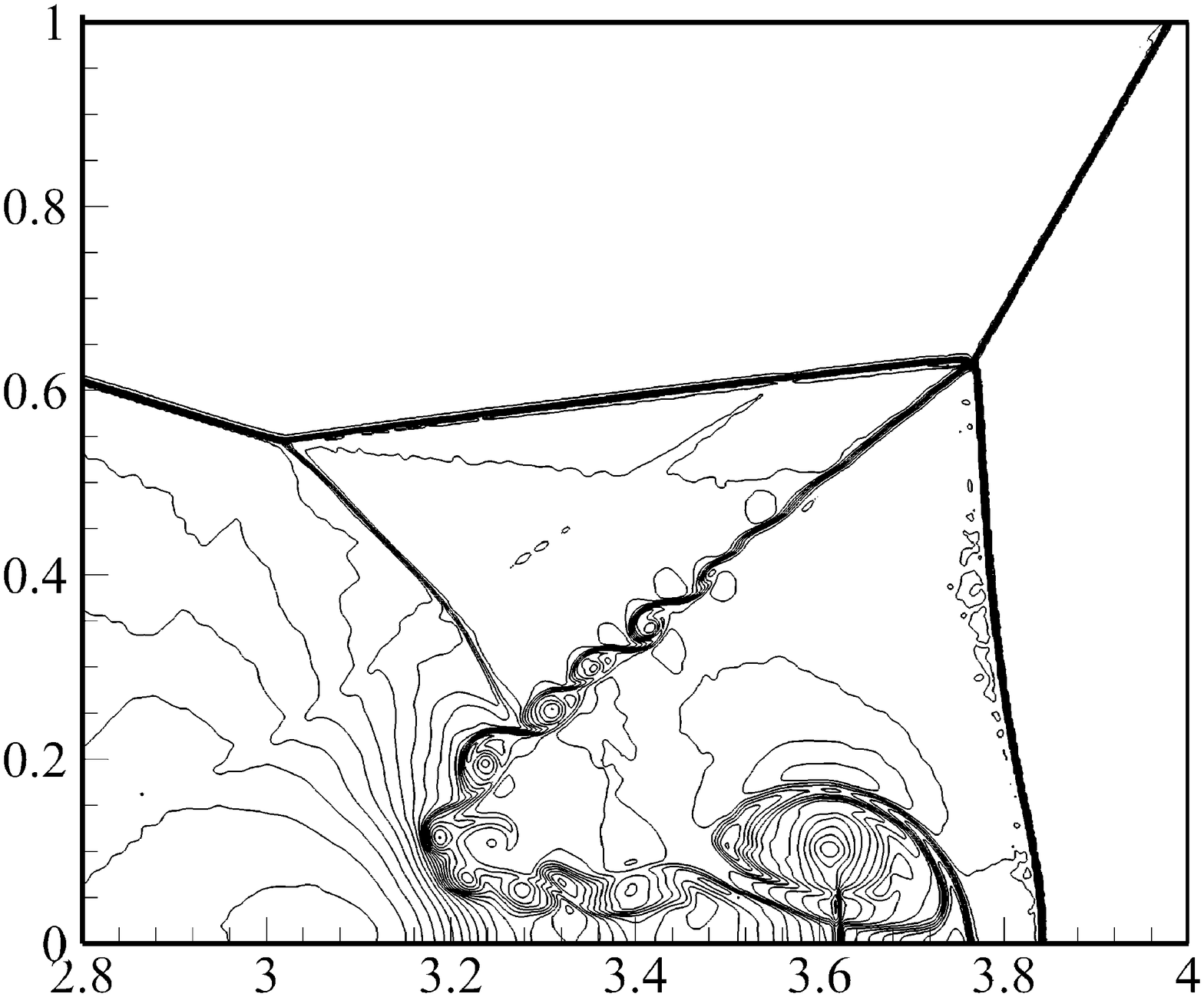}}
  \subfigure[WENO-Z7]{
  \label{FIG:DMR_Density_WENO_Z7_enlarge}
  \includegraphics[width=5.5 cm]{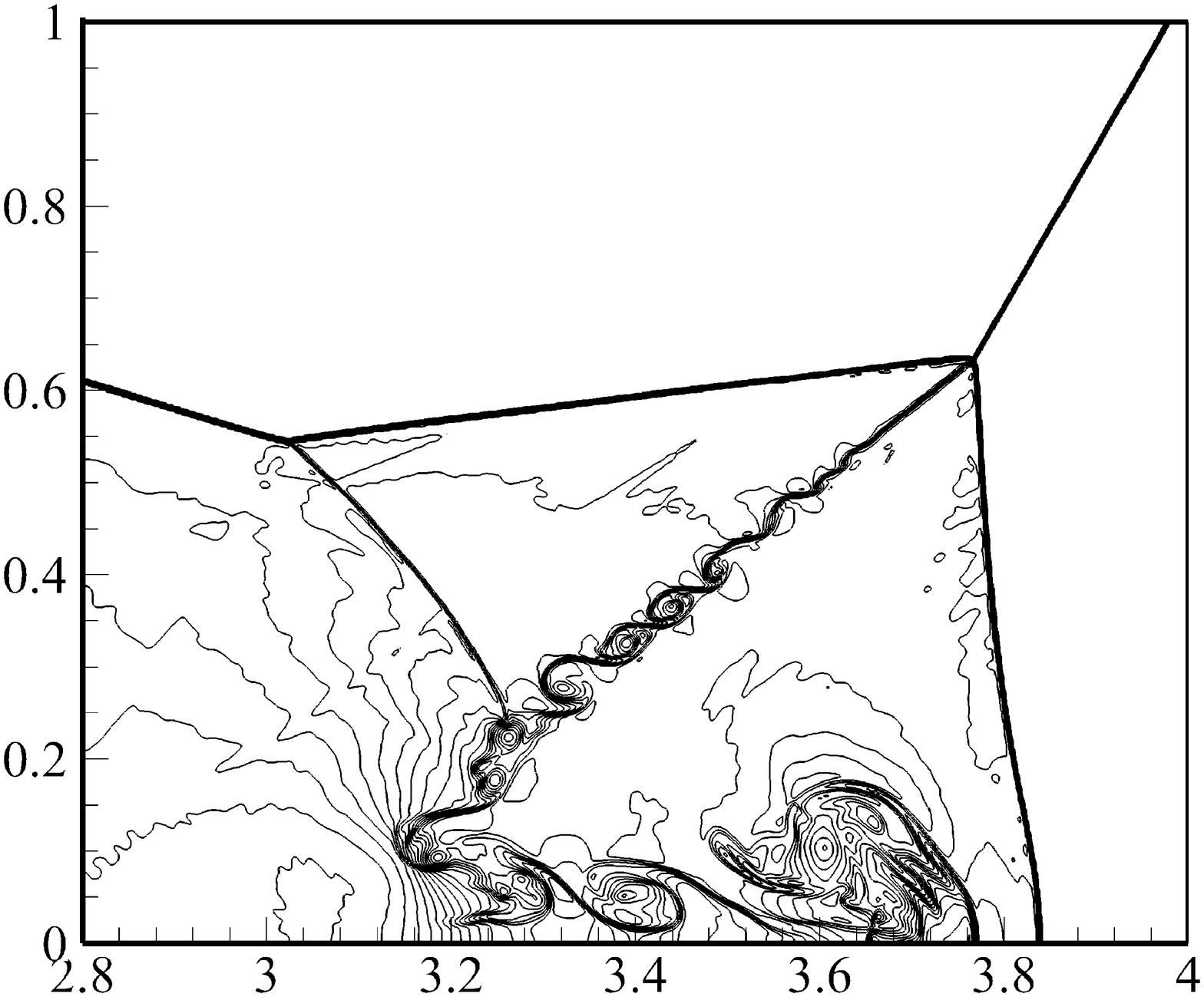}}
  \subfigure[ENO-AO7]{
  \label{FIG:DMR_Density_ENO_AO7_enlarge}
  \includegraphics[width=5.5 cm]{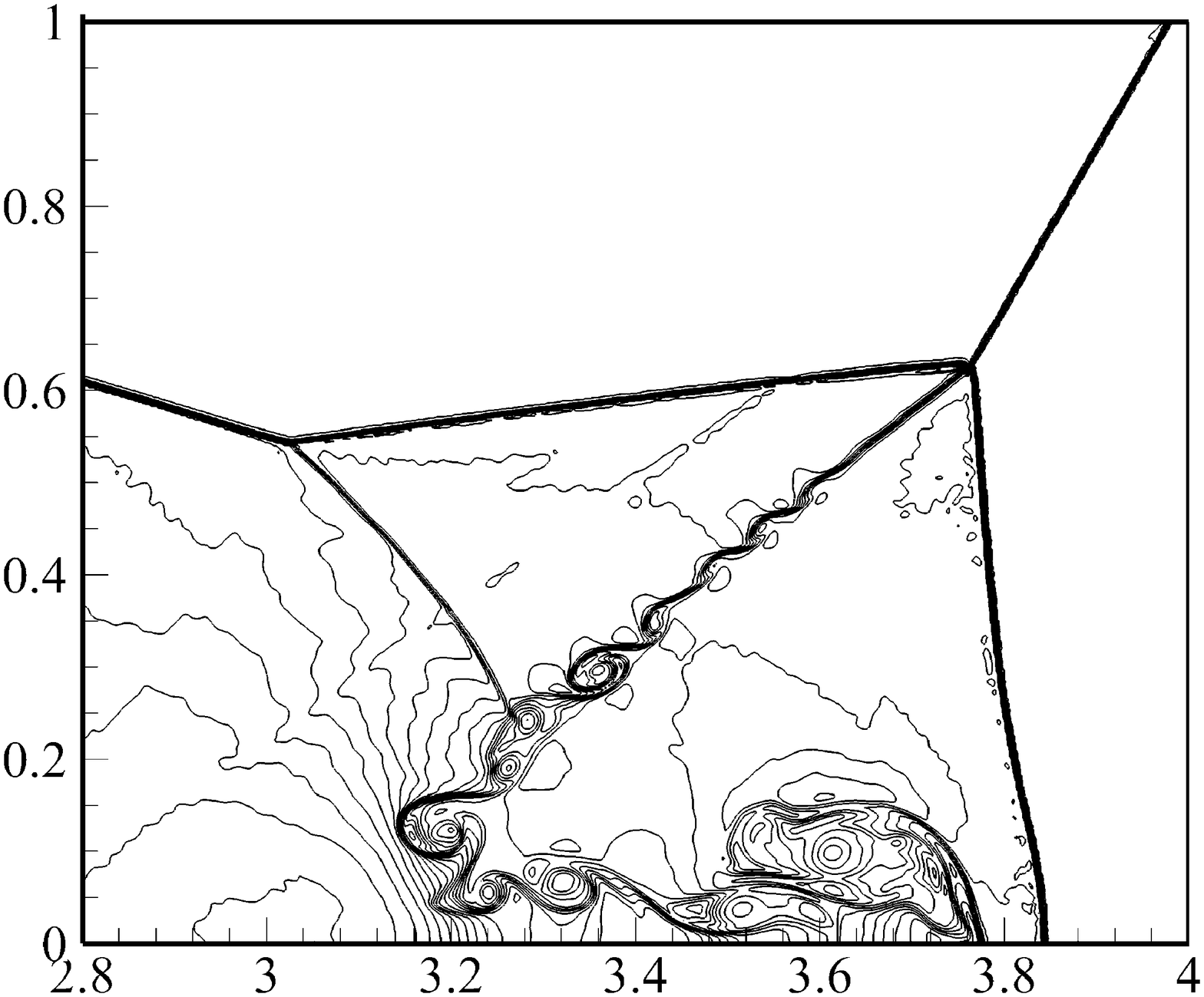}}
  \caption{The zoom-in view of the density contours of the double Mach reflection problem at $t=0.28$ calculated 
  by WENO-Z and ENO-AO with $\Delta x=\Delta y=1/400$.
  The density contours contain 40 equidistant contours from 2 to 22.}
\label{FIG:DMR_enlarge}
\end{figure}

This case was originally proposed by Woodward and Colella \cite{Woodward1984JCP}
and was frequently used to test the shock-capturing and high-fidelity properties of high-order schemes.
The computational domain is $[0, 4]\times[0, 1]$, the ratio of specific heats $\gamma=1.4$,
and the initial condition is given as
\begin{equation*}
  \left(\rho, u, v, p\right)=
    \begin{cases}
      \left(1.4, 0, 0, 1\right) & \text{if } x>\frac{1}{6}+\frac{y}{\sqrt{3}}, \\
      \left(8, 8.25sin(60^\circ), -8.25cos(60^\circ), 116.5\right), & \text{otherwise},
    \end{cases}
\end{equation*}
which describes a Mach 10 oblique shock 
inclined at an angle of $60^\circ$ to the horizontal direction.
On the left, the post-shocked states are imposed;
On the right, nonreflective boundary conditions are implemented;
On the top, the boundary conditions are determined by the exact motion of the oblique shock;
On the bottom, nonreflective boundary conditions are imposed if $x\le\frac{1}{6}$,
otherwise reflective wall boundary conditions are imposed.
Figs. \ref{FIG:DMR} and \ref{FIG:DMR_enlarge} respectively shows 
the entire view and the zoom-in view of the density contours at $t=0.28$ 
calculated by WENO-Z and ENO-AO with $1601\times401$ mesh points.
If we look at the complex flow structures 
near the double Mach reflection zone as shown by Fig. \ref{FIG:DMR_enlarge},
it is fair to say WENO-Z and ENO-AO perform similarly for this case.

\subsubsection{Rayleigh–Taylor instability}
Rayleigh–Taylor instability problem is a canonical physical problem 
that was widely adopted to test the high-fidelity properties of numerical schemes,
see for example the adaptive discontinuous Galerkin schemes \cite{Remacle2003AdaptiveDG}
and the high-order WENO schemes \cite{Shi2003JCP}.

\begin{figure}
  \centering
  \subfigure[WENO-Z5]{
  \label{FIG:DMR_Density_WENO_Z5}
  \includegraphics[height=7 cm]{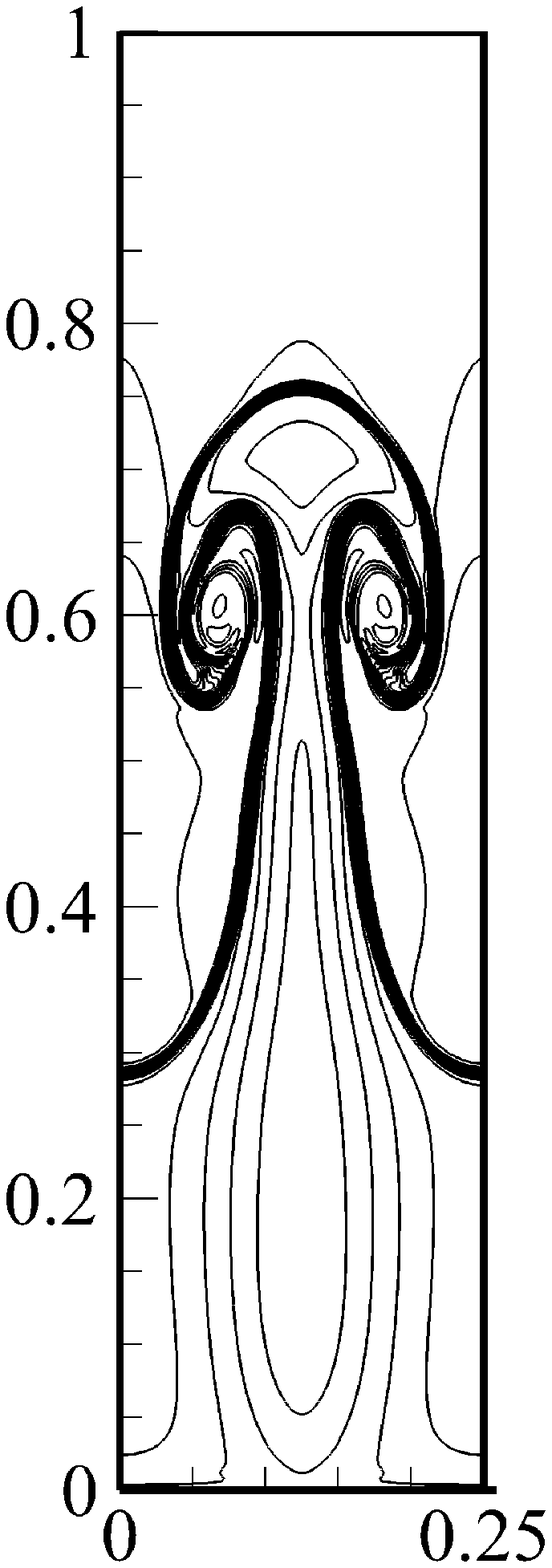}}
  \subfigure[ENO-AO5]{
  \label{FIG:DMR_Density_ENO_AO5}
  \includegraphics[height=7 cm]{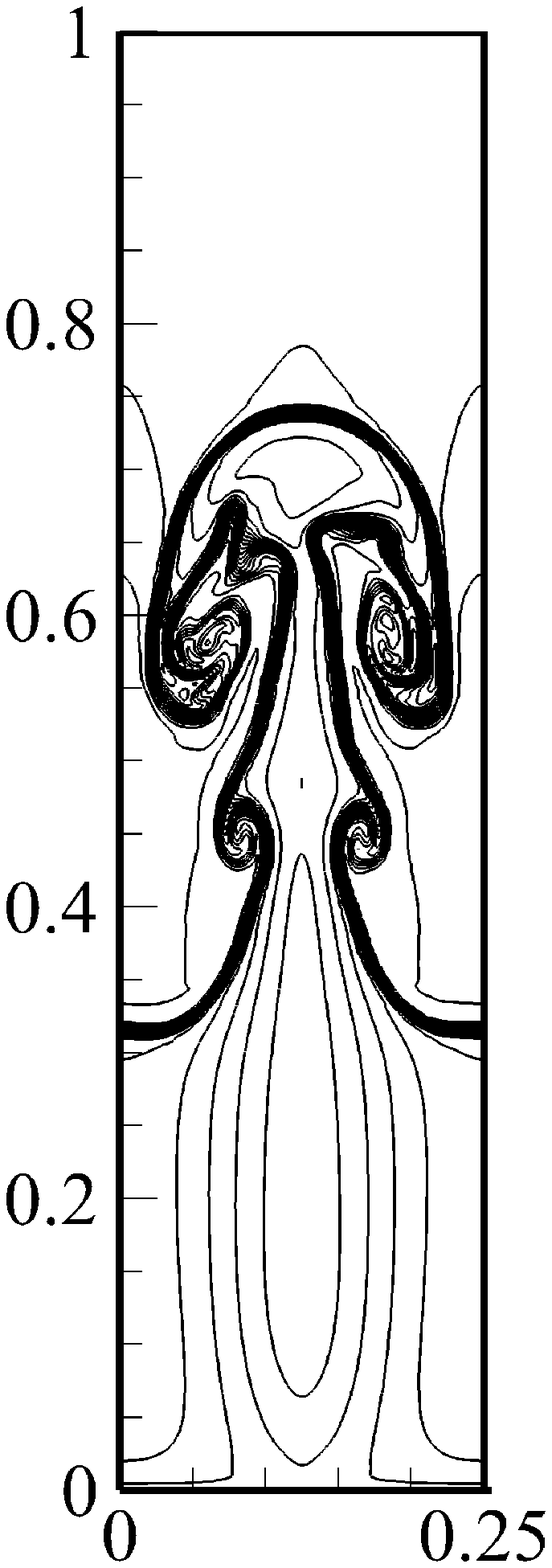}}
  \subfigure[WENO-Z7]{
  \label{FIG:DMR_Density_WENO_Z7}
  \includegraphics[height=7 cm]{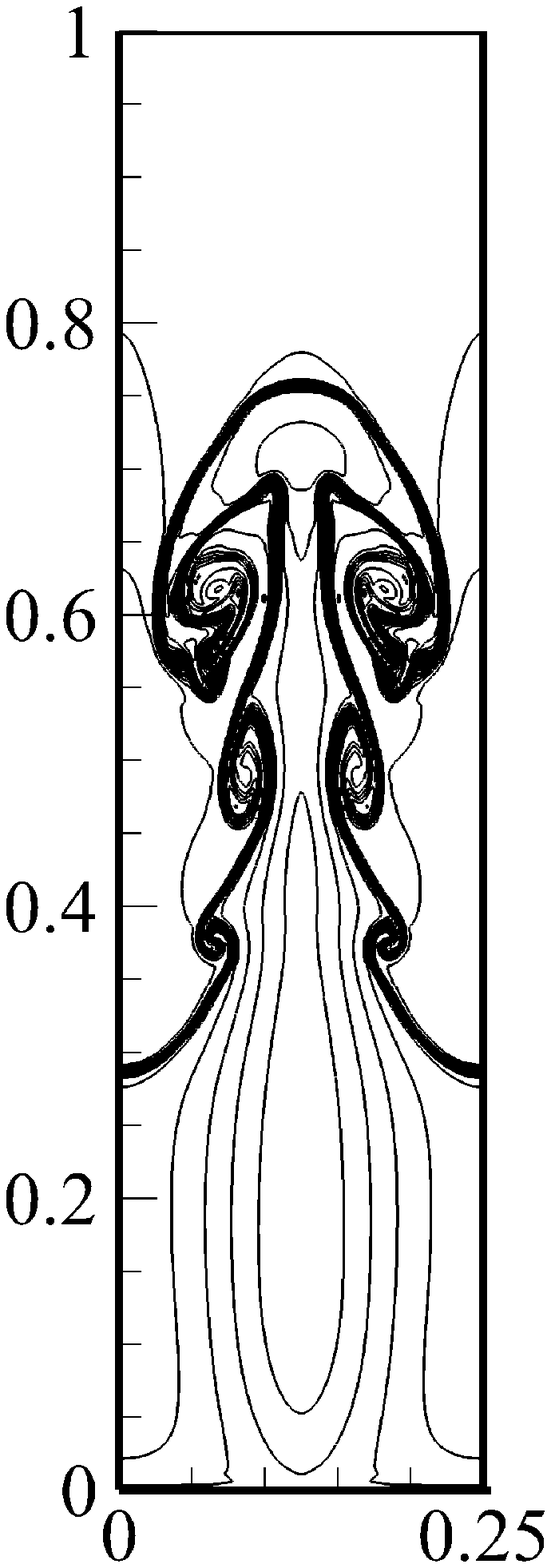}}
  \subfigure[ENO-AO7]{
  \label{FIG:DMR_Density_ENO_AO7}
  \includegraphics[height=7 cm]{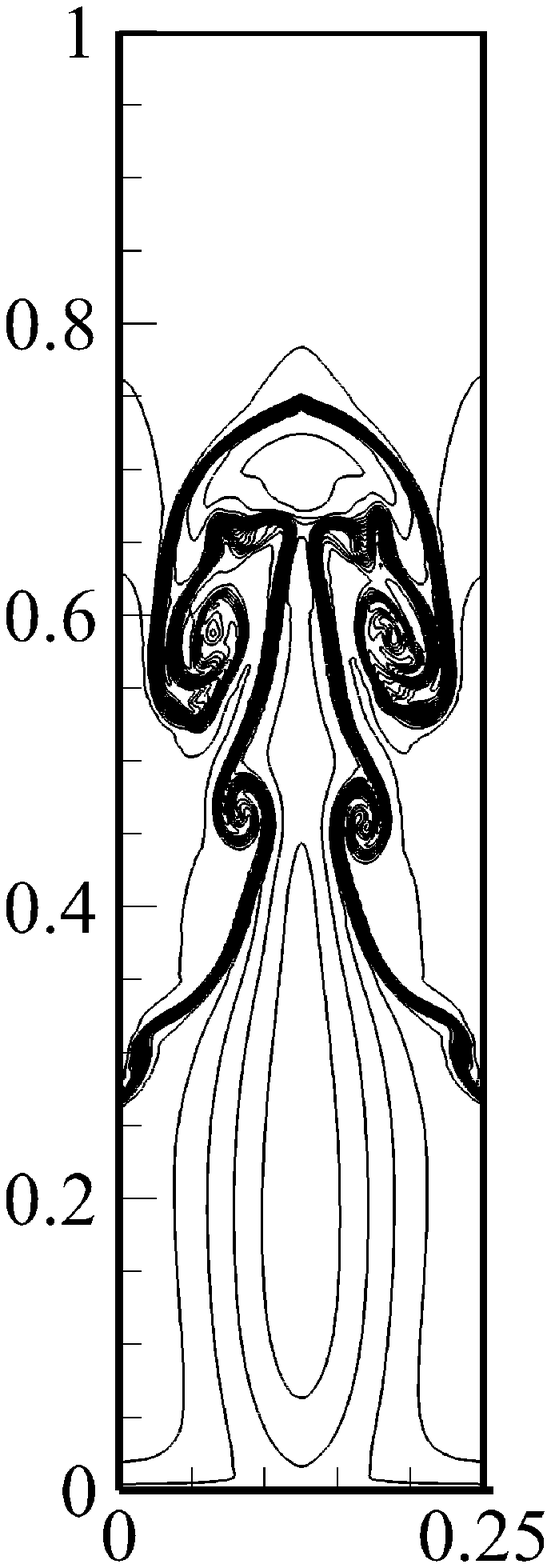}}

  \subfigure[WENO-Z5]{
  \label{FIG:DMR_Density_WENO_Z5}
  \includegraphics[height=7 cm]{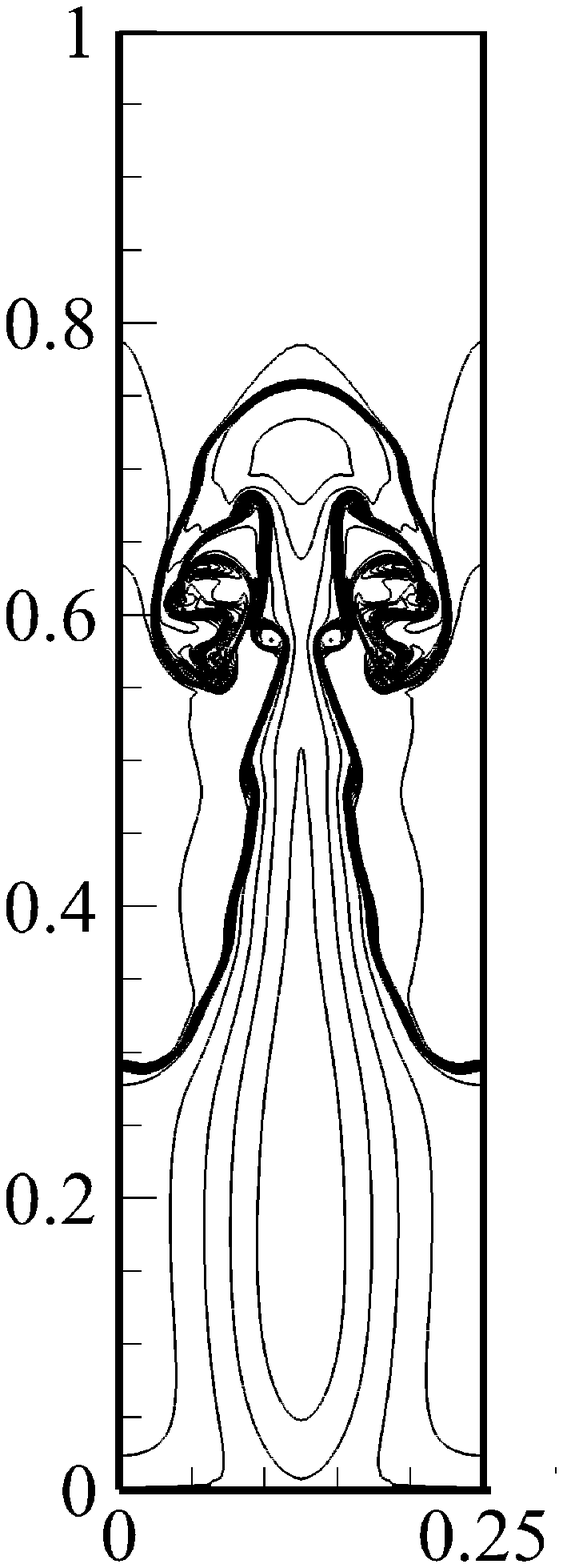}}
  \subfigure[ENO-AO5]{
  \label{FIG:DMR_Density_ENO_AO5}
  \includegraphics[height=7 cm]{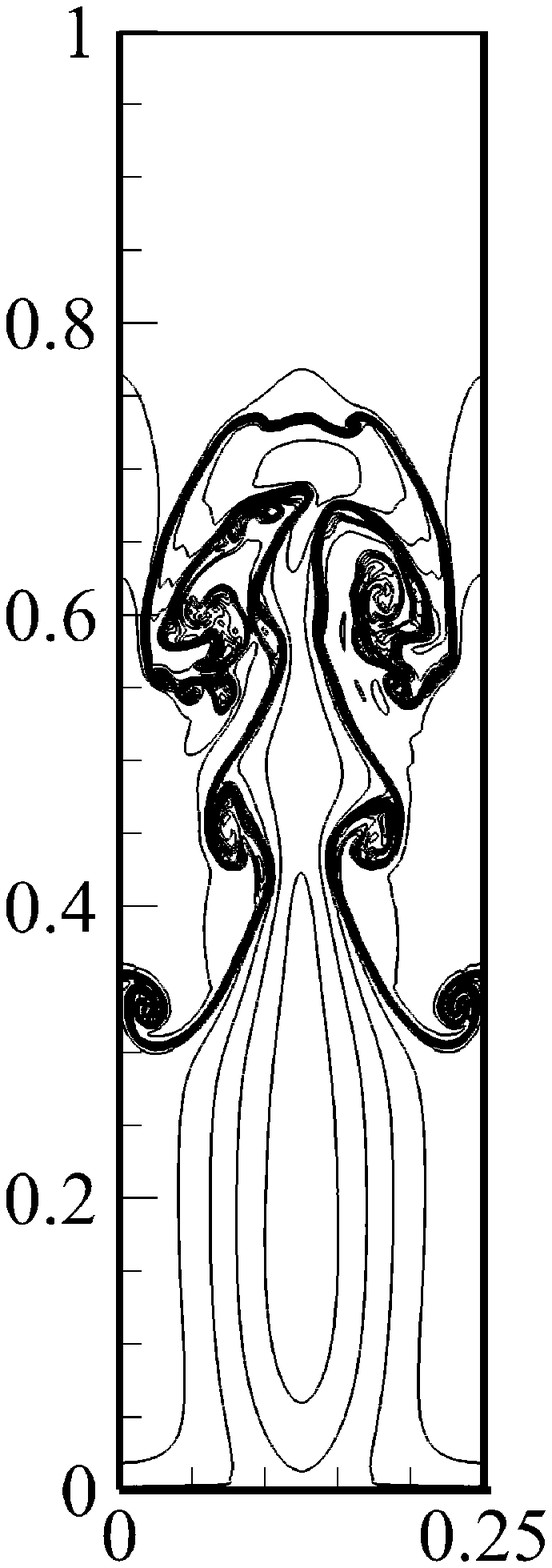}}
  \subfigure[WENO-Z7]{
  \label{FIG:DMR_Density_WENO_Z7}
  \includegraphics[height=7 cm]{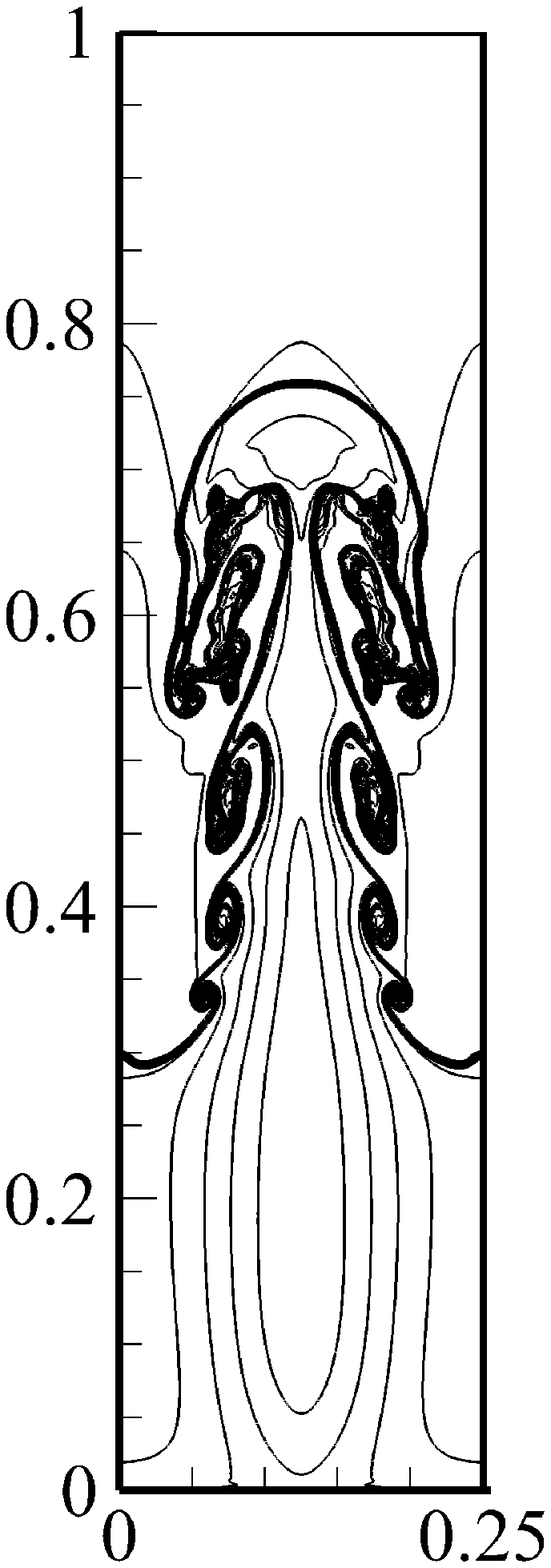}}
  \subfigure[ENO-AO7]{
  \label{FIG:DMR_Density_ENO_AO7}
  \includegraphics[height=7 cm]{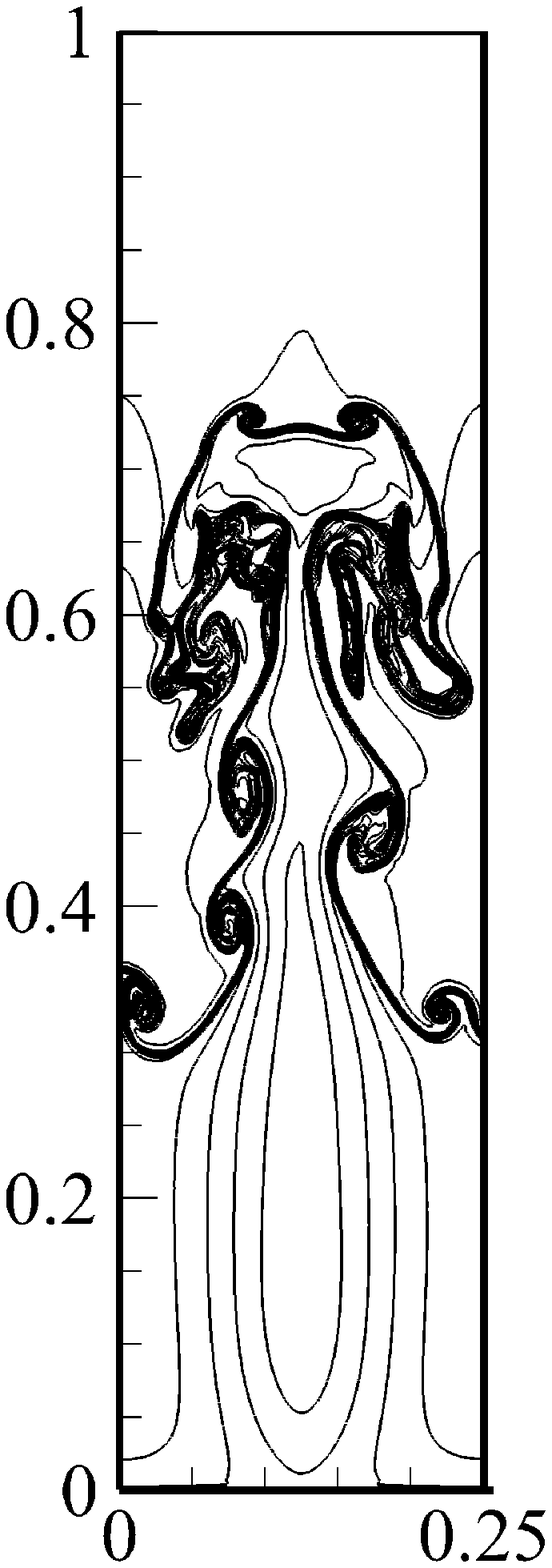}}
  
  \caption{The density contours of the Rayleigh–Taylor instability problem at $t=1.95$ calculated 
  by WENO-Z and ENO-AO with $\Delta x=\Delta y=1/512$ (top row) and $\Delta x=\Delta y=1/1024$ (bottom row).
  The density contours contain 20 equidistant contours from 0.9 to 2.2.}
\label{FIG:RTI}
\end{figure}

We follow the setup of Shi \emph{et al.} \cite{Shi2003JCP}.
The source term $\mathbf{S}=[0,0,\rho,\rho v]^T$ is added to 
the right hand side of the 2D Euler euqations, Eq. (\ref{Eq:2D_Euler_equations}).
The computational domain is $[0,0.25]\times[0,1]$, the ratio of specific heats $\gamma=5/3$,
and the initial condition is given by
\begin{equation*}
  \left(\rho, u, v, p\right)=
    \begin{cases}
      \left(2, 0, -0.025\sqrt{\frac{\gamma p}{\rho}}cos(8\pi x), 2y+1\right) & \text{if } y<0.5, \\
      \left(1, 0, -0.025\sqrt{\frac{\gamma p}{\rho}}cos(8\pi x), y+1.5\right), & \text{otherwise}.
    \end{cases}
\end{equation*}
When the simulations start,
the development of Rayleigh–Taylor instabilities will induce fingering structures which 
cause the mixing of light fluid and heavy fluid.
Fig. \ref{FIG:RTI} shows the density contours at $t=1.95$ calculated 
by WENO-Z and ENO-AO schemes with $\Delta x=\Delta y=1/512$ and $\Delta x=\Delta y=1/1024$ respectively.
ENO-AO5 obviously captures more details of the instabilities than WENO-Z5 does.
It is hard to tell which one of WENO-Z7 and ENO-AO7 performs better, because each one
captures some structures that are missed by another one.
We note that ENO-AO schemes lose symmetry because small truncation errors may cause 
the switch between different stencils.
Many other low-dissipation shock-capturing schemes also encounter symmetry-breaking phenomena
which can be fixed by using a special technique to implemented the computer programs 
\cite{Fleischmann2019symmetry}. 
If we have no special requirements, it is not necessary to enforce the symmetry 
because such instability phenomena are unsymmetric in nature.

\section{Conclusions}
We propose a class of ENO-AO schemes which select candidate stencils 
with unequal sizes by novel smoothness indicators.
The embedded stencils range from first-order all the way up to the designed high-order,
and the new smoothness indicators can adaptively pick up the optimal stencil
according to the smoothness of local solutions.
Therefore, the constructed ENO-AO schemes achieve optimal convergence rates in smooth regions
while maintaining oscillation-free near discontinuities.
The abundant numerical examples show that ENO-AO5 outperforms WENO-Z5 almost for all test cases.
ENO-AO7 and WENO-Z7 are evenly matched for most cases, but ENO-AO7 is more robust for discontinuities 
and performs much better for the Riemann problem that describes the interaction of contact discontinuities.
In summary, the new smoothness indicators provide an efficient way 
to construct high-fidelity shock-capturing schemes with adaptive order.

\section*{Acknowledgement}
H. S. would like to acknowledge the financial support of the National Natural Science Foundation of China (Contract No. 11901602).

\begin{appendix}
\section{The C code of 7th-order ENO-AO7 reconstruction}\label{AppendixA}
The following C code is used to reconstruct $\hat{f}_{j+1/2}^+$ from $f_{j-3}^+,...,f_{j+3}^+$.
In the code, $Fi[k]=f_{j+k-4}^+$.
\small{
\begin{lstlisting}
double ENO_AO7_FluxReconstruction(double Fi[])
{
  int k,Index;
  double DFL,DFR,MinDF,DF[9];
  static double delta=1E-5;
  
  DFL=fabs(-Fi[0]+7*Fi[1]-21*Fi[2]+35*Fi[3]-35*Fi[4]+21*Fi[5]-7*Fi[6]+Fi[7]);
  DFR=fabs(-Fi[1]+7*Fi[2]-21*Fi[3]+35*Fi[4]-35*Fi[5]+21*Fi[6]-7*Fi[7]+Fi[8]);
  DF[8]=MIN(DFL,DFR);
  if(DF[8]<=delta) 
  return (-3*Fi[1]+25*Fi[2]-101*Fi[3]+319*Fi[4]+214*Fi[5]-38*Fi[6]+4*Fi[7])/420;
  
  DFL=fabs(Fi[1]-6*Fi[2]+15*Fi[3]-20*Fi[4]+15*Fi[5]-6*Fi[6]+Fi[7]);
  DFR=fabs(Fi[2]-6*Fi[3]+15*Fi[4]-20*Fi[5]+15*Fi[6]-6*Fi[7]+Fi[8]);
  DF[7]=MIN(DFL,DFR);
  if(DF[7]<=delta) 
  return (Fi[2]-8*Fi[3]+37*Fi[4]+37*Fi[5]-8*Fi[6]+Fi[7])/60.0;
  
  DFL=fabs(Fi[0]-6*Fi[1]+15*Fi[2]-20*Fi[3]+15*Fi[4]-6*Fi[5]+Fi[6]);
  DFR=fabs(Fi[1]-6*Fi[2]+15*Fi[3]-20*Fi[4]+15*Fi[5]-6*Fi[6]+Fi[7]);
  DF[6]=MIN(DFL,DFR);
  if(DF[6]<=delta) 
  return (-Fi[1]+7*Fi[2]-23*Fi[3]+57*Fi[4]+22*Fi[5]-2*Fi[6])/60.0;
  
  DFL=fabs(-Fi[1]+5*Fi[2]-10*Fi[3]+10*Fi[4]-5*Fi[5]+Fi[6]);
  DFR=fabs(-Fi[2]+5*Fi[3]-10*Fi[4]+10*Fi[5]-5*Fi[6]+Fi[7]);
  DF[5]=MIN(DFL,DFR);
  if(DF[5]<=delta) 
  return (2*Fi[2]-13*Fi[3]+47*Fi[4]+27*Fi[5]-3*Fi[6])/60.0;
  
  DFL=fabs(Fi[2]-4*Fi[3]+6*Fi[4]-4*Fi[5]+Fi[6]);
  DFR=fabs(Fi[3]-4*Fi[4]+6*Fi[5]-4*Fi[6]+Fi[7]);
  DF[4]=MIN(DFL,DFR);
  if(DF[4]<=delta) 
  return (-Fi[3]+7*Fi[4]+7*Fi[5]-Fi[6])/12.0;
  
  DFL=fabs(Fi[1]-4*Fi[2]+6*Fi[3]-4*Fi[4]+Fi[5]);
  DFR=fabs(Fi[2]-4*Fi[3]+6*Fi[4]-4*Fi[5]+Fi[6]);
  DF[3]=MIN(DFL,DFR);
  if(DF[3]<=delta) 
  return (Fi[2]-5*Fi[3]+13*Fi[4]+3*Fi[5])/12.0;
  
  DFL=fabs(-Fi[2]+3*Fi[3]-3*Fi[4]+Fi[5]);
  DFR=fabs(-Fi[3]+3*Fi[4]-3*Fi[5]+Fi[6]);
  DF[2]=MIN(DFL,DFR);
  if(DF[2]<=delta) 
  return (-Fi[3]+5*Fi[4]+2*Fi[5])/6.0;
  
  DFL=fabs(Fi[3]-2*Fi[4]+Fi[5]);
  DFR=fabs(Fi[4]-2*Fi[5]+Fi[6]);
  DF[1]=MIN(DFL,DFR);
  if(DF[1]<=delta) 
  return 0.5*(Fi[4]+Fi[5]);
  
  DFL=(fabs(Fi[3]-Fi[4])+fabs(Fi[2]-Fi[3]))/2;
  DFR=(fabs(Fi[5]-Fi[4])+fabs(Fi[6]-Fi[5]))/2;
  DF[0]=MIN(DFL,DFR);
  if(DF[0]<=delta) 
  return Fi[4];
  
  MinDF=DF[0];
  Index=0;
  for(k=1; k<9; k++)
  {
    if(DF[k]<=MinDF)
    {
      MinDF=DF[k];
      Index=k;
    }
  }
  
  switch(Index)
  {
  case 0:
   return Fi[4];
  case 1:
   return 0.5*(Fi[4]+Fi[5]);
  case 2:
   return (-Fi[3]+5*Fi[4]+2*Fi[5])/6.0;
  case 3:
   return (Fi[2]-5*Fi[3]+13*Fi[4]+3*Fi[5])/12.0;
  case 4:
   return (-Fi[3]+7*Fi[4]+7*Fi[5]-Fi[6])/12.0;
  case 5:
   return (2*Fi[2]-13*Fi[3]+47*Fi[4]+27*Fi[5]-3*Fi[6])/60.0;
  case 6:
   return (-Fi[1]+7*Fi[2]-23*Fi[3]+57*Fi[4]+22*Fi[5]-2*Fi[6])/60.0;
  case 7:
   return (Fi[2]-8*Fi[3]+37*Fi[4]+37*Fi[5]-8*Fi[6]+Fi[7])/60.0;
  default:
   return (-3*Fi[1]+25*Fi[2]-101*Fi[3]+319*Fi[4]+214*Fi[5]-38*Fi[6]+4*Fi[7])/420;
  }
}
\end{lstlisting}
}
\end{appendix}
\bibliography{mybibfile}

\end{document}